\documentclass[10pt,leqno]{amsart}

\usepackage{amsmath,amsthm,amsfonts,amscd,amssymb}
\usepackage{hyperref,wasysym}
\usepackage{verbatim}
\usepackage{amssymb}
\usepackage{graphicx}
\usepackage{color,xcolor}

\newtheorem{thm}{Theorem}[section]

\theoremstyle{definition}

\theoremstyle{remark}
\newtheorem{rem}{Remark}[section]

\begin{document}

\title[Structured lattices]{Structured lattices and their applications to security}

\author{Lenny Fukshansky}
\author{Camilla Hollanti}\thanks{This work was supported in part by the Finnish Research Council (Grant \#351271) and in part by the Business Finland Co-Innovation Consortium (Grant \#5845483). R. Y. Njah Nchiwo was supported by the Magnus Ehrnrooth Foundation and the Finnish Academy of Science and Letters, Finland.}
\author{Rahinatou Y. Njah Nchiwo}

\address{Department of Mathematics, 850 Columbia Avenue, Claremont McKenna College, Claremont, CA 91711, USA}
\email{lenny@cmc.edu}
\address{Department of Mathematics and Systems Analysis, Aalto University, P.O. Box 11100, FI-00076 Aalto, Finland}
\email{camilla.hollanti@aalto.fi, rahinatou.njahepousenchiwo@aalto.fi}


\begin{abstract} 
Euclidean lattices are an interesting object of study in many regards and can have a rich structure arising from various constructions, e.g., from number field extensions. A particularly interesting class is the one of well-rounded lattices, as they relate to the well-known densest sphere packing problem in geometry, theta function minimization, and the famous Minkowski and Woods conjectures. In addition to being an important mathematical object in their own right, lattices also play a central role in many applications. This paper offers a survey of structured lattices and discusses their recent applications in lattice-based cryptography and secure wireless communications.  Our goal is to spark the interest of  mathematicians and adjacent communities in these fascinating topics in the intersection of lattices, number theory, cryptography, and wireless communications.  
\end{abstract}

\maketitle

\def\A{{\mathcal A}}
\def\B{{\mathcal B}}
\def\C{{\mathcal C}}
\def\D{{\mathcal D}}
\def\F{{\mathcal F}}
\def\x{{\mathcal H}}
\def\I{{\mathcal I}}
\def\J{{\mathcal J}}
\def\K{{\mathcal K}}
\def\L{{\mathcal L}}
\def\M{{\mathcal M}}
\def\N{{\mathcal N}}
\def\O{{\mathcal O}}
\def\cP{{\mathcal P}}
\def\R{{\mathcal R}}
\def\s{{\mathcal S}}
\def\V{{\mathcal V}}
\def\W{{\mathcal W}}
\def\X{{\mathcal X}}
\def\Y{{\mathcal Y}}
\def\H{{\mathcal H}}
\def\Z{{\mathcal Z}}
\def\OO{{\mathcal O}}
\def\BB{{\mathbb B}}
\def\cee{{\mathbb C}}
\def\EE{{\mathbb E}}
\def\FF{{\mathbb F}}
\def\Nn{{\mathbb N}}
\def\pee{{\mathbb P}}
\def\que{{\mathbb Q}}
\def\real{{\mathbb R}}
\def\zed{{\mathbb Z}}
\def\hyp{{\mathbb H}}
\def\aa{{\mathfrak a}}
\def\HH{{\mathfrak H}}
\def\TT{{\mathfrak \Theta}}
\def\qbar{{\overline{\mathbb Q}}}
\def\eps{{\varepsilon}}
\def\ahat{{\hat \alpha}}
\def\bhat{{\hat \beta}}
\def\gt{{\tilde \gamma}}
\def\h{{\tfrac12}}
\def\be{{\boldsymbol e}}
\def\bei{{\boldsymbol e_i}}
\def\bff{{\boldsymbol f}}
\def\ba{{\boldsymbol a}}
\def\bb{{\boldsymbol b}}
\def\bc{{\boldsymbol c}}
\def\bm{{\boldsymbol m}}
\def\bn{{\boldsymbol n}}
\def\bk{{\boldsymbol k}}
\def\bi{{\boldsymbol i}}
\def\bl{{\boldsymbol l}}
\def\bq{{\boldsymbol q}}
\def\bu{{\boldsymbol u}}
\def\bt{{\boldsymbol t}}
\def\bs{{\boldsymbol s}}
\def\bv{{\boldsymbol v}}
\def\bw{{\boldsymbol w}}
\def\bx{{\boldsymbol x}}
\def\bX{{\boldsymbol X}}
\def\bz{{\boldsymbol z}}
\def\bwy{{\boldsymbol y}}
\def\bY{{\boldsymbol Y}}
\def\bL{{\boldsymbol L}}
\def\baa{{\boldsymbol\alpha}}
\def\bbb{{\boldsymbol\beta}}
\def\bgg{{\boldsymbol\gamma}}
\def\bet{{\boldsymbol\eta}}
\def\bxi{{\boldsymbol\xi}}
\def\bo{{\boldsymbol 0}}
\def\bol{{\boldkey 1}_L}
\def\ep{\varepsilon}
\def\p{\boldsymbol\varphi}
\def\q{\boldsymbol\psi}
\def\rank{\operatorname{rank}}
\def\aut{\operatorname{Aut}}
\def\lcm{\operatorname{lcm}}
\def\sgn{\operatorname{sgn}}
\def\spn{\operatorname{span}}
\def\md{\operatorname{mod}}
\def\Norm{\operatorname{Norm}}
\def\dim{\operatorname{dim}}
\def\det{\operatorname{det}}
\def\Vol{\operatorname{Vol}}
\def\rk{\operatorname{rk}}
\def\Gal{\operatorname{Gal}}
\def\Sym{\operatorname{Sym}}
\def\WR{\operatorname{WR}}
\def\WO{\operatorname{WO}}
\def\GL{\operatorname{GL}}
\def\SL{\operatorname{SL}}
\def\pr{\operatorname{pr}}
\def\Tr{\operatorname{Tr}}
\def\dd{\partial}
\def\itt{\operatorname{int}}
\def\Ar{\operatorname{Area}}
\def\Aut{\operatorname{Aut}}

\tableofcontents

\section{Lattices}\label{sec:lattice}

\subsection{Introduction to lattice theory}
\label{intro}

The theory of Euclidean lattices has its origins in the work of Lagrange and Gauss, who studied lattices in the context of Kepler's sphere packing conjecture and the arithmetic of quadratic forms. Major advances in the theory came with Minkowski's development of geometry of numbers and further connections to number theory, convex and discrete geometry, algebraic geometry, optimization, geometric combinatorics and many other fields of mathematics. Today, lattice theory enjoys a central place within mathematics and its applications. Some of the major breakthroughs of the recent decades included T. Hales's \& S. Ferguson's proof of Kepler's conjecture~\cite{hales}, O. Musin's proof of the kissing number conjecture in dimension~$4$~\cite{musin}, and M. Viazovska et al. proof of the optimal sphere packing conjecture in dimensions~$8$ and~$24$, \cite{viaz1} and~\cite{viaz2}. The goal of this survey paper is to give an overview of theory of Euclidean lattices and some of its recent developments with a view towards applications in coding theory and cryptography. We will especially focus on properties and constructions of important classes of lattices with additional structure (\emph{e.g.}, well-rounded, eutactic, perfect), which play a central role in optimization problems and applications. For comprehensive sources on the theory of Euclidean lattices, we refer the reader to the classical books by Conway \& Sloane~\cite{conway}, Gruber \& Lekkerkerker~\cite{lek}, and Martinet~\cite{martinet}. The first two authors of this paper have recently edited a special collection of research articles on ``Euclidean lattices: theory and applications" (Communications in Mathematics, vol 31, no 2, 2023) which is surveyed in~\cite{lf_ch_survey}.

Throughout this paper we view $\real^n$ as a Euclidean space with respect to the usual Euclidean inner-product $\left<\ ,\ \right>$ and the corresponding Euclidean norm~$\|\bx\| := \sqrt{\left< \bx, \bx \right>}$ for every $\bx \in \real^n$. A {\it lattice} $L \subset \real^n$ of {\it rank} $r \leq n$ (denoted $\rk(L) = r$) is a discrete subgroup of $\real^n$ which is co-compact in the $r$-dimensional subspace $V(L) := \spn_{\real} L$. This is equivalent to saying that there exists a collection of $\real$-linearly independent vectors $\ba_1,\dots,\ba_r \in L$ such that
$$L = \spn_{\zed} \{ \ba_1,\dots,\ba_r \} = \left\{ \sum_{i=1}^r c_i \ba_i : c_1,\dots,c_r \in \zed \right\}.$$
The collection $\ba_1,\dots,\ba_r$ is a basis for $L$ and we refer to the $n \times r$ matrix $A = (\ba_1 \dots \ba_r)$ as a {\it basis matrix} for $L$, \emph{i.e.}, $L = A\zed^r$. For any $U \in \GL_r(\zed)$, $AU$ is another basis matrix for $L$. The subspace $V(L) \subseteq \real^n$ can be identified with~$\real^r$, so from here on we will talk about lattices of {\it full rank} in $\real^n$, meaning that $r=n$. The space $\L_n$ of full-rank lattices in $\real^n$ can be identified with $\GL_n(\real) \backslash \GL_n(\zed)$, the set of orbits of $\GL_n(\real)$ under the action of $\GL_n(\zed)$ by right multiplication.

Two lattices $L_1,L_2 \subset \real^n$ are said to be {\it similar}, denoted $L_1 \sim L_2$, if there exists $\alpha \in \real^+$, the group of positive real numbers, and $U \in \O_n(\real)$, the $n \times n$ real orthogonal group, such that $L_2 = \alpha U L_1$. This is an equivalence relation, and the space $\s_n$ of similarity classes is $(\real^+ \times \O_n(\real)) / \L_n$, the set of orbits of the space of full-rank lattices under the action of the group $\real^+ \times \O_n(\real)$ by left multiplication. A lattice $L$ is called {\it integral} if for every $\bx,\bwy \in L$, $\left< \bx, \bwy \right> \in \zed$ and a lattice is called {\it arithmetic} if it is similar to an integral lattice.

The group of isometries of a full-rank lattice $L \subset \real^n$ is $\O(L) := \{ U \in \O_n(\real) : UL = L \}$, which is compact as a subset of $\GL_n(\real)$ with respect to the Euclidean metric topology. The automorphisms of $L$ are isometries given by integer linear transformations, so the {\it automorphism group} of $L$ is $\Aut(L) := \GL_n(\real) \cap \O(L)$. Since $\Aut(L)$ is the intersection of a discrete subgroup with a compact set, it is finite. In fact, in all but seven exceptional dimensions the lattice with the largest (with respect to size) automorphism group is~$\zed^n$, which is the signed permutation group consisting of $2^n n!$ elements: it is generated by independent coordinate sign changes and permutations. When $n=2,4,6,7,8,9,10$ more symmetric lattices with even larger automorphism groups exist. For example, $|\Aut(\zed^2)|=8$ whereas the hexagonal lattice
$$\Lambda_h = \begin{pmatrix} 1 & 1/2 \\ 0 & \sqrt{3}/2 \end{pmatrix} \zed^2$$
has $12$ automorphisms. This being said, most lattices have just two automorphisms: multiplication by $\pm 1$. The automorphism group is an invariant of the similarity class of a lattice.

The {\it determinant} (also called {\it co-volume}) of a lattice $L = A\zed^n$ is defined as $\det(L) := \sqrt{\det(A^{\top} A)}$ and is equal to the volume the quotient group $\real^n/L$. This is an invariant of $L$ which does not depend on the choice of a basis matrix $A$. 
 In other words, $\det(L)$ is the volume of any {\it fundamental domain}, \emph{i.e.} a measurable full set of coset representatives of $L$ in $\real^n$. One example of a fundamental domain is a fundamental parallelotope $\{ c_1 \ba_1 + \dots + c_n \ba_n : 0 \leq c_i < 1\ \forall\ 1 \leq i \leq n \}$, corresponding to the basis matrix $A = (\ba_1 \dots \ba_n)$. Another object more intrinsically dependent on $L$ is its {\it Voronoi cell}
 $$\V(L) := \{ \bx \in \real^n : \|\bx\| \leq \|\bx-\bwy\|\ \forall\ \bwy \in L \}.$$
While $\V(L)$ is not a fundamental domain of $L$, it is the closure of a fundamental domain, and hence its volume is still equal to $\det(L)$. Now, $\real^n = \bigcup_{\bx \in L} (\bx + \V(L))$, where two distinct translates $\bx_1 + \V(L)$ and $\bx_2 + \V(L)$ can intersect only at the boundary. A lattice $L$ is called {\it unimodular} if $\det(L) =1$.
 
 For a given full-rank lattice $L = A\zed^n \subset \real^n$, its {\it dual lattice} is defined to be
 $$L^* := \{ \bx \in \real^n : \left< \bx,\bwy \right> \in \zed\ \forall\ \bwy \in L \} = (A^{-1})^{\top} \zed^n,$$
 then $\det(L^*) = 1/\det(L)$. If $L$ is integral then $L \subseteq L^*$, hence unimodular integral lattices are {\it self-dual}, \emph{i.e.}, $L = L^*$.
 
We also define the {\it successive minima} $0 < \lambda_1(L) \leq \dots \leq \lambda_n(L)$ of a full-rank lattice $L$ in $\real^n$ to be
$$\lambda_i(L) := \min \left\{ t \in \real^+ : \dim_{\real} \left( \spn_{\real} ( L \cap \BB_n(t) \right) \geq i \right\},$$
where $\BB_n(t)$ is a ball of radius $t$ centered at the origin in $\real^n$. In particular, the first successive minimum is the norm of a shortest nonzero vector in $L$. Additionally, the {\it covering radius} (also called the {\it inhomogeneous minimum}) of $L$ is defined as
$$\mu(L) := \min \left\{ t \in \real^+ : L + \BB_n(t) = \real^n \right\}.$$
The celebrated Minkowski's Successive Minima Theorem (see, \emph{e.g.}, Chapter~2, \S~9 of~\cite{lek}) gives bounds on the product of successive minima of $L$:
\begin{equation}
\label{mink_s}
\frac{2^n \det(L)}{n! \omega_n} \leq \prod_{i=1}^n \lambda_i(L) \leq \frac{2^n \det(L)}{\omega_n},
\end{equation}
where $\omega_n$ is the volume of a unit ball in $\real^n$. An immediate implication of~\eqref{mink_s} is Minkowski Convex Body Theorem (see, \emph{e.g.}, Chapter~2, \S~5 of~\cite{lek}):
\begin{equation}
\label{mink_c}
\lambda_1(L) \leq 2 \left( \frac{\det(L)}{\omega_n} \right)^{1/n}.
\end{equation}
We define the set of {\it minimal vectors} of $L$ to be $S(L) := \{ \bx \in L : \|\bx\| = \lambda_1(L) \}$.

The {\it sphere packing} associated to $L$ is constructed by inscribing a ball of radius $\lambda_1(L)/2$ into each translate of the Voronoi cell $\V(L)$. The {\it density} of this packing is the proportion of the space occupied by the spheres, which is the same as ratio of the volume of the ball and volume of the Voronoi cell into which it is inscribed; it can be computed as
$$\delta(L) := \frac{\omega_n \lambda_1(L)^n}{2^n \det(L)} \leq 1.$$
This is a continuous function on the space of lattices $\L_n$ which is constant on any given similarity class, hence we can think of it as a continuous function of the space $\s_n$ of similarity classes. The objective of the {\it lattice packing problem} in a given dimension is to find a lattice that maximizes $\delta(L)$. Solutions to the lattice packing problem are only known in dimensions $1 \leq n \leq 9$ (with $n=9$ case being very recent still unpublished work by Dutour Sikiri\'c and van Woerden) and $n=24$ (see Chapter~1 of~\cite{conway}). The celebrated Minkowski-Hlawka theorem (see, \emph{e.g.}, Chapter~1 of~\cite{conway}) asserts that in every dimension $n \geq 2$ there exists a full-rank lattice $L \subset \real^n$ such that
$$\delta(L) \geq \frac{\zeta(n)}{2^{n-1}},$$
where $\zeta(n)$ the value of the Riemann zeta-function at $n$. The proof of this theorem, however, is not constructive, and in all but finitely many dimensions (up to $1000$ or so) no constructions of lattices satisfying the Minkowski-Hlawka bound are known. This being said, the lower bound $\zeta(n)/2^{n-1}$ in the Minkowski-Hlawka theorem can be improved. The most significant improvement is due to B.~Klartag~\cite{klartag}, who established (also non-constructively) the lower bound $cn^2/2^{n}$ for some universal constant $c > 0$ in 2025. Notice that maximizing lattice packing density in dimension $n \geq 2$ is equivalent to determining the value of the Hermite's constant
$$\gamma_n := \max_{L \subset \real^n} \frac{\lambda_1(L)}{\det(L)^{1/n}}.$$
While the value of $\gamma_n$ is only known in the few dimensions mentioned above, an upper bound on it follows, for instance, from Minkowski's theorem~\eqref{mink_c}.

A linearly independent collection of vectors $\ba_1,\dots,\ba_n \in L$ is said to {\it correspond to successive minima} if $\|\ba_i\| = \lambda_i$ for each $1 \leq i \leq n$. Finding the successive minima is equivalent to the \emph{shortest independent vector problem} mentioned in Section \ref{sec:lbc}, which is known to be NP-hard.

Such a collection is not unique, but there are only finitely many of them in a given lattice. These vectors are known to form a basis for $L$ in dimensions $n \leq 3$, however in dimensions $n \geq 4$ they do not necessarily form a basis. Consider, for example the lattice
$$L_1 = \spn_{\zed} \left\{ \be_1, \be_2,\be_3, \frac{1}{2} \sum_{i=1}^4 \be_i \right\} \subset \real^4,$$
where $\be_i$ are the standard basis vectors. Then
$$\left\{ \be_1,\be_2,\be_3,\frac{1}{2} \sum_{i=1}^4 \be_i \right\} \text{ and } \left\{ \be_1,\be_2,\be_3,\be_4 \right\}$$
are both collections of vectors in $L_1$ corresponding to successive minima, however the first one forms a basis whereas the second does not. Similarly, the lattice
\begin{equation}
\label{L2}
L_2 = \spn_{\zed} \left\{ \be_1, \be_2,\be_3, \be_4, \frac{1}{2} \sum_{i=1}^5 \be_i \right\} \subset \real^5
\end{equation}
does not have a collection of vectors corresponding to successive minima that would form a basis. These observations raise a natural question: what is the shortest basis in $L$? Hermite inequality (see, \emph{e.g.} Theorem~2.2.1 of~\cite{martinet}) guarantees that a lattice $L$ of rank $n$ has a basis $\ba_1,\dots,\ba_n$ such that $\lambda_1(L) = \|\ba_1\| \leq |\ba_2\| \leq \dots \leq \|\ba_n\|$ and
\begin{equation}
\label{hermite}
\prod_{i=1}^n \|\ba_i\| \leq \left( \frac{4}{3} \right)^{\frac{n(n-1)}{2}} \det (L).
\end{equation}
On the other hand, Hadamard inequality (see, \emph{e.g.} Theorem~2.1.1 of~\cite{martinet}) states that the {\it orthogonality defect} of this basis is
\begin{equation}
\label{hadamard}
\nu(\ba_1,\dots,\ba_n) := \frac{\prod_{i=1}^n \|\ba_i\|}{\det(L)} \geq 1.
\end{equation}
Indeed, this basis is orthogonal if and only if $\nu(\ba_1,\dots,\ba_n) = 1$ and
$$\frac{\lambda_1(L)^n}{\det(L)} \leq \nu(\ba_1,\dots,\ba_n) \leq \left( \frac{4}{3} \right)^{\frac{n(n-1)}{2}},$$
implying that maximizing $\frac{\lambda_1(L)^n}{\det(L)}$ to achieve $\gamma_n^n$ entails maximizing the orthogonality defect.

We also want to mention two related optimization problems on lattices: the {\it lattice covering problem} and the {\it kissing number problem}. The covering configuration associated to the lattice $L$ is constructed by circumscribing a sphere of radius $\mu(L)$ around each translate of the Voronoi cell, and the {\it thickness} of this covering is the ratio of the volume of the ball and volume of the Voronoi cell around which it is circumscribed; it can be computed as
$$\TT(L) := \frac{\omega_n \mu(L)^n}{\det(L)} \geq 1.$$
Again, this is a continuous function of the space $\s_n$ of similarity classes of lattices. The objective of the {\it lattice covering problem} in a given dimension is to find a lattice that minimizes $\TT(L)$. Solutions to the lattice covering problem are only known in dimensions $1 \leq n \leq 5$ (see Chapter~1 of~\cite{conway}). Finally, the kissing number problem on lattices asks for the maximal number of spheres centered at points of a lattice $L$ in $\real^n$ that can touch the sphere centered at $\bo$. This is equivalent to asking for a lattice with maximal number of minimal vectors in a given dimension, \emph{i.e.}, the kissing number of $L$ is $|S(L)|$. The answer is known in dimensions $1 \leq n \leq 9$ and $n=24$ (see Chapter~1 of~\cite{conway}).
\bigskip

\subsection{Well-rounded and related classes}
\label{WR}

A lattice $L \subset \real^n$ is called {\it well-rounded} (abbreviated WR) if $\lambda_1(L) = \dots = \lambda_n(L)$, which is equivalent to saying that 
\begin{equation}
\label{WR_def}
\real^n = \spn_{\real} S(L).
\end{equation}
It is important to remark that WR condition~\eqref{WR_def} is not equivalent to the condition $L = \spn_{\zed} S(L)$, as demonstrated by the example $L_2$ in~\eqref{L2} above: this second condition is strictly stronger than~\eqref{WR_def} for $n \geq 5$; if it holds, we say that $L$ is {\it generated by minimal vectors} (for $n \leq 4$, all WR lattices are generated by their minimal vectors). Further, for all $n \geq 10$ it is possible for a lattice $L$ to be generated by minimal vectors while not containing a basis of minimal vectors: this was first demonstrated for $n \geq 11$ by Conway and Sloane~\cite{sloane} and then extended to $n \geq 10$ by Martinet and Sch\"urmann~\cite{bases-3}. On the other hand, for all $n \leq 9$ lattices generated by minimal vectors contain a basis of minimal vectors (see \cite{bases-1}, \cite{bases-2}, \cite{bases-3}).

WR property is preserved under similarity, hence we can speak of WR similarity classes, of which there are infinitely many in each dimension $n \geq 2$. WR lattices appear in many different contexts in number theory, geometry, combinatorics and optimization. In particular, the space of WR lattices forms a ``spine" ($\SL_n(\zed)$-equivariant deformation retract) for the space of all lattices in $\real^n$, which is useful for cohomology computations~\cite{ash1}, \cite{pettet}, \cite{lizhen_ji}. Further, WR lattices appear prominently in regard to Minkowski conjecture~\cite{mcmullen}. Let $L$ be a unimodular lattice. For each $\bx \in L$, define the multiplicative norm $N(\bx) = |x_1 \cdots x_n|$. Notice that $N$ is preserved under the left-multiplication action by the diagonal group
$$\A_n := \left\{ \begin{pmatrix} a_1 & \dots & 0 \\ \vdots & \ddots & \vdots \\ 0 & \dots & a_n \end{pmatrix} : a_i > 0, \prod_{i=1}^n a_i = 1 \right\},$$
in other words, $N(\bx) = N(A\bx)$ for any $A \in \A_n$. Minkowski conjectured that for any unimodular lattice $L \subset \real^n$,
$$\sup_{\bx \in \real^n} \inf_{\bwy \in L} N(\bx-\bwy) \leq \frac{1}{2^n}.$$
This conjecture was originally motivated by the study of certain ``approximation properties" of algebraic integers in number fields (see~\cite{bayer_nebe})). Minkowski proved his conjecture for $n=2$; up until 2005, the conjecture was proved in dimensions $n \leq 5$. In his seminal paper~\cite{mcmullen}, C. McMullen proved this conjecture for $n=6$ (see~\cite{mcmullen} for references to the earlier work). He follows the Remak-Davenport approach (see Section 27.1 of~\cite{gruber} for details and history), splitting the Minkowski's conjecture into two statements from which it follows:

$(W_n)$ For any lattice $L \subset \real^n$, there exists $A \in \A_n$ such that $AL$ is WR.

$(C_n)$ For any unimodular WR lattice $L \subset \real^n$, $\mu(L) \leq \mu(\zed^n) = \sqrt{n}/2$.

\noindent
In the direction of $(W_n)$, McMullen established that if the closure of the orbit of $L$ under the action of $\A_n$ is bounded, then it contains a WR lattice. This, along with $(C_n)$, turns out to be enough to establish Minkowski's conjecture. The second part $(C_n)$ is known as A. C. Woods' covering conjecture~\cite{woods}; prior to McMullen's work, it has been proved in dimensions $n \leq 6$ and has since been proved in all dimensions $n \leq 10$ (see~\cite{raka} and references within). On the other hand, Woods' conjecture has been disproved in dimensions $n \geq 30$ by Regev, Shapira and Weiss~\cite{weiss-1},  whose work has been extended by Chen and Xu to show that the conjecture fails for $n \geq 24$~\cite{chen_xu}. This, however, does not necessarily mean that Minkowski conjecture in those dimensions is not true.

Another important class of lattices (related to WR lattices at least in spirit) is semi-stable lattices: $L \subset \real^n$ is called {\it semi-stable} if for every sublattice $M \subseteq L$,
$$\det(M)^{1/\rk(M)} \leq \det(L)^{1/\rk(L)}.$$
Semi-stable lattices were first introduced in the context of reduction theory, where this condition was taken to heuristically suggest that the successive minima of $L$ are not too far apart (see~\cite{andre} and the excellent survey paper of Casselman~\cite{casselman} on semi-stable lattices, which in particular provides many references and the history of development of this subject). This, however, does not mean that WR lattices are necessarily semi-stable: this statement is only true for $n=2$, whereas for $n \geq 3$ the sets of semi-stable and WR lattices are independent (see, \emph{e.g.},~\cite{lf_st} for explicit examples of non-stable WR lattices in $\real^3$), although they do have an intersection. Similarly to the WR lattices, semi-stable lattices are also well-distributed among the orbits of the diagonal group action on the space of lattices. Specifically, in~\cite{weiss-2} the authors showed that, analogously to McMullen's observation about WR lattices, if the closure of the orbit of $L$ under the action of $\A_n$ is bounded, then it contains a semi-stable lattice. Remarkably, in~\cite{solan} Solan strengthened this observations for both, WR and semi-stable lattices, proving that for any lattice $L \subset \real^n$, there exist $A, B \in \A_n$ such that $AL$ is WR and $BL$ is semi-stable. In particular, this establishes the conjecture $(W_n)$ from above.

\begin{figure}[t]
\centering
\includegraphics[scale=0.4]{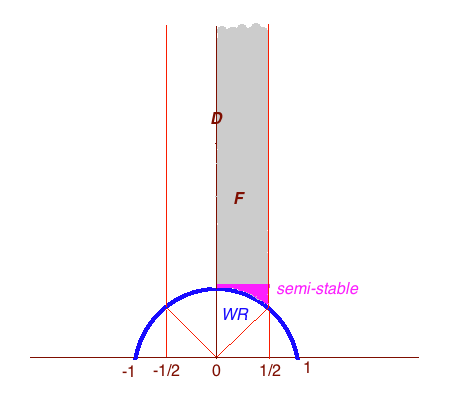}
\caption{Similarity classes of lattices in $\real^2$ with WR and semi-stable subregions marked by colors.}\label{fig:domain}
\end{figure}

The two-dimensional distribution of WR and semi-stable lattices can be described very explicitly and deserves some attention due to its connection to the parameterization of elliptic curves; our brief exposition follows \cite{pasha1}, \cite{pasha2}, \cite{pasha3}. Let $\hyp = \{ \tau = a+bi : b \geq 0 \} \subset \cee$ be the upper half-plane, and let
$$\D := \{ \tau = a + bi \in \hyp : -1/2 < a \leq 1/2, |\tau| \geq 1 \}.$$
Let
$$\F := \{ \tau = a + bi \in \hyp : 0 \leq a \leq 1/2, |\tau| \geq 1 \},$$
so, loosely speaking, $\F$ is ``half" of $\D$. Every point $\tau = a + bi \in \F$ can be identified with a lattice
$$\Gamma_{\tau} := \begin{pmatrix} 1 & a \\ 0 & b \end{pmatrix} \zed^2$$
in $\real^2$. Every planar lattice $L$ is similar to a unique lattice of the form $\Gamma_{\tau}$ for some $\tau \in \F$, hence we can say that the similarity class of $L$ is represented by $\tau$. Thus, $\F$ can be thought of as the space of similarity classes of lattices in~$\real^2$ (see Figure~\ref{fig:domain}). WR similarity classes correspond to the circular arc $\{ \tau \in \F : |\tau|=1 \}$ and semi-stable similarity classes correspond to the set $\{ \tau = a+bi \in \F : b \leq 1 \}$. On the other hand, the full domain $\D$ can be viewed as the space of isomorphism classes of elliptic curves: a point $\tau$ corresponds to the isomorphism class of the elliptic curve given by the complex torus $\cee/\Gamma_{\tau}$, where we are identifying $\cee$ with $\real^2$ and thinking of $\Gamma_{\tau}$ as $\spn_{\zed} \{1,\tau\} \subset \cee$. This being said, while the lattices $\Gamma_{\tau}$ and $\Gamma_{\bar{\tau}}$ are similar, the corresponding elliptic curves are not isomorphic: instead, the two elliptic curves have conjugate $j$-invariants, since $j(-\bar{\tau}) = \overline{j(\tau)}$  (here $j$ is Klein's modular $j$-function).

Further, the question of distribution of WR sublattices of a given lattice $L \subset \real^2$ has been investigated by several authors via analysis of the properties of the so-called WR zeta-function
$$\zeta_{\WR}(s) = \sum_{k=1}^{\infty} a_k k^{-s},$$
where $a_k$ is the number of WR sublattices of $L$ of index $k$ and $s$ is a complex variable. Information about the position and order of the pole, as well as the residue at the pole of this function can be used along with Wiener-Ikehara Tauberian theorem and its later variations to establish the order of growth of the counting function $\sum_{k \leq n} a_k$, the number of WR sublattices of $L$ of index at most $n$ as $n \to \infty$. The result depends on whether the lattice $L$ is arithmetic or not. Specifically, combining the results of~\cite{lf_wr} and \cite{stefan}, we obtain:
$$\sum_{k \leq n} a_k = \left\{ \begin{array}{ll}
O(n \log n) & \mbox{if $L$ is arithmetic,} \\
O(n) & \mbox{if $L$ is not arithmetic,}
\end{array}
\right.$$
as $n \to \infty$. More detailed asymptotic results on the summatory function $\sum_{k \leq n} a_k$ have later been obtained in~\cite{baake}.

There are several other contexts in which WR lattices have been investigated, \emph{e.g.} in connection with the Frobenius problem~\cite{frob}. Most importantly, WR lattices are key in discrete optimization and applications. In particular, the lattice packing problem can be restricted to WR lattices without loss of generality, \emph{i.e.}, any solution to the lattice packing problem in every dimension $n \geq 2$ has to be a WR lattice. In fact, more is true. A lattice $L \subset \real^n$ with $m = |S(L)|$ is called {\it eutactic} if there exist positive real numbers $c_1,\dots,c_m$ such that
\begin{equation}
\label{eut_crit}
\|\bx\|^2 = \sum_{i=1}^m c_i \left< \bx,\bwy_i \right>^2,
\end{equation}
where $S(L) = \{ \bwy_1,\dots,\bwy_m \}$. The lattice $L$ is {\it strongly eutactic} if $c_1 = \dots = c_m$. On the other hand, if the space of $n \times n$ real symmetric matrices $\Sym_n(\real)$ can be represented as
$$\Sym_n(\real) = \spn_{\real} \left\{ \bwy_i \bwy_i^{\top} : \bwy_i \in S(L) \right\},$$
then the lattice $L$ is called {\it perfect}. Both, perfect and eutactic properties are preserved under similarity. While the sets of eutactic and perfect lattices are independent (there are perfect non-eutactic and eutactc non-perfect lattices), both of them are WR, but not necessarily semi-stable: an example of a perfect non-stable lattice in dimension $8$ has been obtained by Y. Kim~\cite{ykim}, where he also proved that all other perfect lattices in dimension $\leq 8$ are semi-stable. A lattice is called {\it extreme} if it is a local maximum of the packing density function in its dimensions, and a classical theorem of Voronoi (see, \emph{e.g.}, Theorem~3.4.6 of~\cite{martinet}) states that a lattice is extreme if and only if it is perfect and eutactic. It is well known that perfect lattices are necessarily arithmetic, hence so are extreme lattices. In every dimension $n \geq 2$ there are only finitely many eutactic and finitely many perfect similarity classes. For instance, up to similarity in $\real^2$, there are only two eutactic lattices ($\zed^2$ and $\Lambda_h$) and only one perfect ($\Lambda_h$). Further, Bacher proved~\cite{bacher} that the number $p_n$ of perfect similarity classes of lattices in $\real^n$ satisfies the following inequalities for any $\eps > 0$:
$$e^{n^{1-\eps}} < p_n < e^{n^{3+\eps}}.$$
The upper bound has recently been improved by van Woerden~\cite{woerden} to $e^{O(n^2 \log n)}$. The exact value of $p_n$ has so far been published in all dimensions $n \leq 8$ with the $8$-dimensional case being an extensive computational project by Dutour Sikiri\'c, Sch\"urmann and Vallentin~\cite{achill-1} building on the previous results (see Section~6.6 of~\cite{martinet}): they showed that there are $10916$ perfect lattices in~$\real^8$, but only $2408$ of them are eutactic (hence, extreme) by the work of Riener~\cite{riener}. The computational project by Dutour Sikiri\'c and van Woerden to count and classify perfect lattices in $\real^9$ has been ongoing for some years and has recently completed: their work is currently being prepared for publication.

We close this section with the notion of {\it generic well-rounded (GWR)} lattices, which has been considered in~\cite{LevinGWR,camilla}: a full-rank WR lattice $L \subset \real^n$ is called GWR if $|S(L)|=2n$. One immediate example of a GWR lattice is $\zed^n$. GWR lattices can have a relatively high packing density while maintaining a relatively small kissing number, which is a useful property for some coding theory applications we detail below. One steady source of GWR lattices are the nearly orthogonal lattices. Given an ordered basis $B = \{ \bb_1,\dots,\bb_n \}$ for a lattice $L$, define a sequence of angles $\theta_1,\dots,\theta_{n-1}$ so that each $\theta_i$ is the angle between $\bb_{i+1}$ and the subspace 
$$\spn_{\real} \{ \bb_1,\dots,\bb_i \}.$$
Then each $\theta_i \in [0,\pi/2]$ and $B$ is called a {\it weakly nearly orthogonal} basis if $\theta_i \geq \pi/3$ for each $1 \leq i \leq n-1$. A basis $B$ is called {\it nearly orthogonal} if every ordering of it is weakly nearly orthogonal. If $L$ has such a basis, we say that $L$ is a nearly orthogonal lattice. Nearly orthogonal lattices were introduced in~\cite{baraniuk}, where they were applied to the problem of image compression. Nearly orthogonal lattices that are also WR (and often GWR) have been studied in~\cite{lf_dk-1}.

\bigskip

\subsection{Algebraic constructions}
\label{constr}

Due to their importance in a variety of theoretical contexts and applications, explicit constructions of WR families of lattices (often with additional properties) are of great interest. In particular, a great deal of attention was devoted to the study of constructions coming from different algebraic settings. To this end, let us start with a number field $K$ of degree $n \geq 2$ over $\que$ and let us write $\Delta_K$ for its discriminant and $\O_K$ for its ring of integers. Assume that $K$ has $r_1$ real embeddings $\sigma_1,\dots,\sigma_{r_1} : K \to \real$ and $r_2$ pairs of complex conjugate embeddings $\tau_1,\bar{\tau}_1,\dots,\tau_{r_2},\bar{\tau}_{r_2} : K \to \cee$. Then $n = r_1 + 2r_2$ and we can define the {\it Minkowski embedding} of $K$ into $\real^n$ by
$$\Sigma_K := \left( \sigma_1,\dots,\sigma_{r_1},\Re(\tau_1),\Im(\tau_1),\dots,\Re(\tau_{r_2}),\Im(\tau_{r_2}) \right) : K \to \real^n.$$
Let $J \subset K$ be a fractional ideal, then $\Sigma_K(J) \subset \real^n$ is a lattice of full rank. Lattices like this are called {\it ideal lattices} via Minkowski embedding. Notice that for any $\alpha,\beta \in K$,
$$\left< \Sigma_K(\alpha), \Sigma_K(\beta) \right> = \Tr_K(\alpha \bar{\beta}),$$
where $\Tr_K$ stands for the number field trace on $K$. Hence, the Euclidean lattice structure on $\Sigma_K(J)$ is induced by the trace of $K$. The theory of ideal lattices in a more general form has been developed by Bayer-Fluckiger, among other authors; see~\cite{bayer}, \cite{bayer_nebe} for a detailed survey of this area.

WR ideal lattices have first been studied in~\cite{lf_kp}, where it was in particular proved that for totally real and totally imaginary number fields, $\Sigma_K(\O_K)$ is WR if and only if $K$ is cyclotomic (see also~\cite{alves}). More generally, let us say that an ideal $J \subseteq \O_K$ is WR if the corresponding ideal lattice $\Sigma_K(J)$ is WR. Infinite families of real and imaginary quadratic number fields containing WR ideals have been constructed in~\cite{lf_kp}. Further, in~\cite{lf_wr-2} it has been proved that for squarefree positive integer $D$, quadratic fields $K(\sqrt{\pm D})$ contain WR ideals when $D$ has a divisor $d$ satisfying
$$\sqrt{\frac{D}{3}} \leq d < \sqrt{D}.$$
This condition is if and only if in the case $K=\que(\sqrt{-D})$. On the other hand, \cite{anitha} establishes that $\que(\sqrt{D})$ contains WR ideals if and only $D$ has a divisor $d$ satisfying
$$\sqrt{\frac{D}{3}} \leq d < \sqrt{3D}.$$
Thus, relatively few of ideal lattices from quadratic number fields are WR. On the other hand, as follows from the results of McMullen~\cite{mcmullen} and Solan~\cite{solan}, any ideal lattice can be ``twisted" into a WR one by the action of the diagonal group $\A_2$. Damir and Karpuk in~\cite{damir_karpuk} investigated properties of the specific bases of ideals that result in the minimal basis of such corresponding WR twist. Further, situations when the canonical basis of an ideal in a quadratic number field can be so twisted have been studied in~\cite{damir-3}. For higher degree number fields, WR ideals have been proved to exist in cyclic cubic and (some) cyclic quartic number fields~\cite{ha-2}. Semi-stable ideal lattices have also been investigated; in particular, infinite families of semi-stable ideal lattices from any real quadratic number field were constructed in~\cite{lf_st} (where it was also proved that a positive proportion of ideal lattices from real quadratic fields are semi-stable) and semi-stable twists of canonical bases of ideals in all quadratic number fields were studied in~\cite{damir-3}. Well-rounded ideal lattices from totally definite quaternion algebras were recently studied in \cite{oggier_WRquat}.

There is a different algebraic construction of lattices actively used in cryptography, which has also received the name of ideal lattices. Let $f(x) \in \zed[x]$ be a monic polynomial of degree $n$ and consider the quotient ring $R(f) := \zed[x]/\left< f(x) \right>$. Define the {\it coefficient embedding} 
$$\rho_f : R(f) \to \zed^n,$$
given by $\rho_f( a_0 + a_1x + \dots + a_{n-1} x^{n-1}) = (a_0,a_1,\dots,a_{n-1})^{\top}$. This is a linear map between two free $\zed$-modules and an ideal $J \subseteq R(f)$ is mapped onto a sublattice $\rho_f(J) \subseteq \zed^n$. We will call such lattices the {\it coefficient ideal lattices}; they generalize the original construction of {\it cyclic lattices} (introduced by Micciancio in~\cite{mic}), which are the special case $\rho_f(J)$ for $f(x) = x^n-1$. These lattices were introduced and studied in the cryptographic context by Lyubashevsky and Micciancio~\cite{lub_mic}. WR cyclic lattices have been investigated in \cite{xsun}, \cite{xsun-err}, \cite{lf_dk}. 

If $f(x)$ is an irreducible polynomial, then the corresponding map $\rho_f$ is a linear isomorphism of $\zed$-modules $R(f)$ and $\zed^n$ and the coefficient ideal lattices $\rho_f(J)$ have full rank (Lemma~3.2 of~\cite{lub_mic}). In the special case when $f(x)$ is $n$-th cyclotomic polynomial of degree $\varphi(n)$, the ring $R(f)$ can be identified with the ring of integers~$\zed[\theta_n]$ of the cyclotomic field $\que(\theta_n)$ via the canonical isomorphism $x \mapsto \theta_n$, where $\theta_n = e^{2\pi i/n}$ is $n$-th primitive root of unity. One can then ask about the relation between the two embeddings, \emph{i.e.}, between the ideal lattice $\Sigma_{\que(\theta_n)}(J)$ and the coefficient ideal lattice $\rho_f(J)$ for a given ideal $J$ in this ring. The linear transformation between the two has been worked out by Batson in~\cite{batson} as we describe in Eq. \eqref{trans} .

The construction of lattices via Minkowski embedding can start from any free $\zed$-module contained in the number field $K$, not only from an ideal in $\O_K$. Let, for instance, $\M(\alpha) = \spn_{\zed} \{ \alpha_1,\dots,\alpha_n \}$, where $\alpha=\alpha_1,\alpha_2,\dots,\alpha_n$ are algebraic conjugates contained in $K$. A special case of this construction when $K$ is a cyclic number field of odd prime degree has been considered in \cite{robson-1}, \cite{robson-2}, where families of WR such lattices (in fact, even having bases of minimal vectors) have been obtained. More general such constructions of GWR nearly orthogonal lattices with bases of minimal vectors and large automorphism groups -- in particular, coming from Pisot numbers -- are currently being explored in~\cite{lf_ek}. More generally, the so-called {\it module lattices} in $\real^{nd}$, $n = [K:\que]$, $d \geq 1$, coming from modules $\M \subset K^d$ via Minkowski embedding $\Sigma_K: K^d \rightarrow \real^{nd}$ were used in \cite{langlois2015worst} in the construction of lattice-based crypto-schemes; see Section~\ref{sec:lbc} for more details.

Another explicit algebraic construction of WR lattices comes from curves over finite fields. Let $F$ be an algebraic function field of a single variable with the finite field $\FF_q$ as its field of constants. Let $\cP = \{ P_0,P_1,\ldots, P_{n-1}\}$ be the set of rational places of $F$ and for each $P_i$, let $v_i$ denote the corresponding normalized discrete valuation. Let  $\O_{\cP}^*$ be the abelian group of all nonzero functions $f\in F$ whose divisor has support contained in the set $\cP$. Then $\sum_{i=0}^{n-1} v_i(f) = 0$ for each $f \in \O_{\cP}^*$. Let $n= |\cP|$, the number of rational places of $F$, and define the homomorphism $\phi_{\cP} : \O_{\cP}^* \to \zed^n$ by 
$$\phi_{\cP}(f) = (v_0(f),v_1(f), \ldots, v_{n-1}(f)).$$ 
Then $L_{\cP} := \phi_\cP(\O_{\cP}^*)$ is a finite-index sublattice of the root lattice
$$A_{n-1} = \left\{ \bx \in \zed^n : \sum_{i=0}^{n-1} x_i = 0 \right\}.$$
These lattices, called {\it function field lattices}, are described in detail in the well known book~\cite{tv} by Tsfasman and Vladut (Chapter 5.4); they were originally introduced in~\cite{rt} by Rosenbloom and Tsfasman, who used this construction to produce asymptotically good families of lattices from the standpoint of packing density. A more systematic investigation of the geometric properties of these lattices was carried out more recently by several authors. In particular, the case of elliptic curves (algebraic function fields of genus $1$) has been considered in~\cite{lf_hm} and~\cite{sha}, where it has been proved that for $n \geq 5$ the lattice $L_{\cP}$ is WR and has a basis of minimal vectors. This, however, is a special case of the more general result of~\cite{bottcher1} on lattices from finite abelian groups discussed below. For fields of higher genus, the case of Hermitian function fields $\FF_q(x,y)$ with $y^q+y=x^{q+1}$ for a prime power $q$ is considered in~\cite{bottcher3}, where it is proved that the corresponding lattice $L_{\cP}$ is WR and generated by minimal vectors. On the other hand, examples of hyperelliptic function fields giving rise to non-WR lattice $L_{\cP}$ are presented in~\cite{la_hs}.

We now turn to explicit algebraic constructions of families of extreme lattices. The most standard of these are the {\it irreducible root lattices} $A_n$, $D_n$, $E_6$, $E_7$, $E_8$ and on some occasions their duals (a lattice is called {\it irreducible} if it is not an orthogonal sum of proper sublattices). We refer the reader to the excellent detailed exposition of the theory of root lattices (as well as related to them Coxeter lattices) in Martinet's book~\cite{martinet} and focus instead on a more recent lesser-known construction. Let $G = \{ 0_G, z_1,\dots,z_n \}$ be an additive abelian group and define
$$L_G := \left\{ \ba \in \zed^n : \sum_{i=1}^n a_i = 0,\ \sum_{i=1}^n a_i z_i = 0_G \right\},$$
which is a lattice of rank $n-1$. B\"ottcher et al.~\cite{bottcher1} proved that the lattice $L_G$ is WR and, in fact, has a basis of minimal vectors for $G \neq \zed/4\zed$. Further, B\"ottcher et al.~\cite{bottcher2} established that $L_G$ is strongly eutactic for all $G$ of odd order or elementary abelian 2-groups. Additionally, Ladisch~\cite{ladisch} proved that $L_G$ is eutactic for all $G \neq \zed/4\zed$. On the other hand, Bacher~\cite{bacher-0} proved that for all $G$ of order at least $9$, the lattice $L_G$ is perfect. Putting these results together, we see that $L_G$ is extreme for all finite abelian groups $G$ of order at least $9$. Additional constructions of strongly eutactic lattices from tight (equiangular) frames and distance transitive graphs have been given in~\cite{deanna1} and~\cite{deanna2}, respectively. We choose, however, not to detail them here since they are somewhat more analytic in nature.

Finally, we mention the notion of \emph{tame} lattices that were introduced in \cite{damir} and motivated by the behavior of the trace pairing over tame cyclic number fields as well as by their ability to serve a method to explicitly construct WR lattices. Tame lattices have a \emph{Lagrangian basis} \cite{conner-lagrangian} and they have Gram matrices of a nice specific form. Tame lattices are known to exist for any tame number field with a prime conductor \cite{wilmar}. GWR lattices arising from tame ones were constructed and studied in \cite{camilla}, including a discussion on their applicability to wireless security, which we discuss in Section~\ref{sec:IT-sec}.
\bigskip

\subsection{Spherical designs and Epstein zeta-function}
\label{e-zeta}

We also briefly discuss spherical designs and their connection to (highly structured) Euclidean lattices. Let $\s_{n-1}(r) \subset \real^n$ be the sphere of radius $r$ centered at the origin in $\real^n$. A finite collection of points $\{ \bx_1,\dots,\bx_m \} \subset \s_{n-1}(r)$ is called a {\it spherical $t$-design} for an integer $t \geq 1$ if for any polynomial $p(y_1,\dots,y_n) \in \real[y_1,\dots,y_n]$ of degree $\leq t$,
$$\frac{1}{m} \sum_{i=1}^m p(\bx_i) = \int_{\s_{n-1}(r)} p(\bwy) d\bwy,$$
where the measure is normalized so that $\int_{\s_{n-1}(r)} d\bwy = 1$. Spherical $t$-designs have been originally introduced by Delsarte, Goethals and Seidel~\cite{sph_designs}. The existence of spherical $t$-designs in $\real^n$ for any $t,n \geq 1$ was established by Seymour and T. Zaslavsky~\cite{seym_zasl}) with effective bounds on the size of such designs produced by Bondarenko, Radchenko and Viazovska~\cite{bond_rad_viaz}. We state here a convenient criterion for a $\bo$-symmetric set $X \subset \s_{n-1}(r)$ to be a spherical $t$-design for $t=2p$ and $t=2p+1$, $p \geq 1$: such $X$ is a spherical $t$-design if and only if there exists a constant $c_p$ such that 
\begin{equation}
\label{sph_d_crit}
\|\bwy\|^{2p} = \frac{c_p}{r^{2p} |X|} \sum_{\bx \in X} \left< \bwy,\bx \right>^{2p},
\end{equation}
for all $\bwy \in \real^n$. An important class of such symmetric spherical designs comes from sets of minimal vectors in lattices. Indeed, comparing~\eqref{eut_crit} with~\eqref{sph_d_crit}, we see that a lattice $L \subset \real^n$ is strongly eutactic if and only its set of minimal vectors $S(L)$ is a spherical $2$-design.

The connection between spherical designs and lattices has been studied by several authors, starting with the fundamental paper of B. Venkov~\cite{venkov} (see also Chapter 16 of~\cite{martinet} for a nice exposition of the results of~\cite{venkov}). A lattice whose set of minimal vectors is a spherical $4$-design is called {\it strongly perfect}, and Venkov proves that strongly perfect lattices are extreme (hence, perfect, by Voronoi's theorem, making the notation justified). In the same paper, Venkov also produced a list (albeit not a full classification) of strongly perfect lattices in dimensions $n \leq 24$. Various classification results for strongly perfect lattices and lattices carrying even higher-degree spherical designs have appeared since (see, \emph{e.g.},~\cite{nossek} and references within). We refer the reader to the paper~\cite{nebe_venkov} by Nebe for a detailed survey of Venkov's theory of lattices and spherical designs and related contributions by other authors.

Spherical designs play a role in another important optimization problem on lattices. The {\it Epstein zeta-function} of a lattice $L \subset \real^n$ is defined as
$$\Z_L(s) = \sum_{\bx \in L \setminus \{ \bo \}} \|\bx\|^{-2s},$$
for a variable $s \in \cee$. For each lattice $L$ in each dimension $n \geq 2$, this Dirichlet series converges in the half-plane $\Re(s) > n/2$, has a simple pole at $s=n/2$ and admits a meromorphic continuation to the whole complex plane. The classical minimization problem for the Epstein zeta-function considers a fixed real value $s_0 > 0$, $s_0 \neq n/2$, and asks for a lattice $L_0 \subset \real^n$ such that
$$\Z_{L_0}(s_0) = \min \left\{ \Z_L(s_0) : L \subset \real^n \right\}.$$
Besides the intrinsic lattice theory interest, this problem also come up in the work of S. Sobolev~\cite{sobolev} in regards to numerical integration.
In dimension $n=2$, this problem dates back at least to the work of Rankin~\cite{rankin}, Cassels~\cite{cassels_epstein}, Diananda~\cite{diananda} and Ennola~\cite{ennola}, who established that for any such $s_0$ the minimum occurs only at the hexagonal lattice. In dimension $n=3$, the minimization problem was solved by Ennola~\cite{ennola_1} and in dimensions $n= 4,8,24$ by Sarnak and Str\"ombergsson~\cite{sarnak}. It was separately proved by Ry\v{s}kov~\cite{ryskov} that the minimizer of $\Z_L(s)$ as $s \to \infty$ corresponds to the densest lattice packing in $\real^n$.

This last observation makes it especially interesting to look for such minimizers, and that is where spherical designs again make an appearance. Given a lattice $L \subset \real^n$, let us define the {\it spectrum} of $L$ to be the set $\left\{ \|\bx\| : \bx \in L \setminus \{ \bo \} \right\}$. We can write the spectrum as an ordered set of real numbers $\{ 0 < a_1 < a_2 < \dots \}$ and define the $k$-th {\it layer} of $L$ to be
$$\{ \bx \in L : \|\bx\| = a_k \}.$$
With this notation, we can state a remarkable theorem of Coulangeon~\cite{coulangeon}: if every layer of $L$ contains a spherical $4$-design, then $L$ is a minimizer of $\Z_L(s)$ for every real value of $s > n/2$.

The minimization problem for Epstein zeta-function also has an interesting applied connection. A Dirichlet series $F(s) = \sum_{n=1}^{\infty} a_n n^{-s}$ and the corresponding power series $f(z) = \sum_{n=1}^{\infty} a_n z^n$ are connected via the {\it Mellin transform}:
$$\Gamma(s) F(s) = \int_0^{\infty} x^{s-1} f(e^{-x}) dx,$$
where $\Gamma(s)$ is the value of the $\Gamma$-function at $s$. Thus, $\Z_L(s)$ for an integral lattice $L$ corresponds to the power series $\Theta_L(z)$, called the {\it theta-function} of $L$. The corresponding minimization problem for $\Theta_L(z)$ has important implications for maximizing the reliability and security of communications when using lattices to construct coding schemes for wireless communications. We discuss this connection in more detail in Section~\ref{sec:IT-sec}.

\bigskip

\section{Lattice-based cryptography}\label{sec:lbc}
Over the last few decades, interest in lattices has grown due to their applications in cryptography. In this chapter, we discuss lattice-based cryptography (LBC), focusing on paradigms whose security depends on hard mathematical problems involving lattices. We describe several complex problems in LBC and examine their worst-case and average-case hardness. In particular, our main focus will be on  the Learning with Errors (LWE) problem and its variants, and analyze the relationships among them. For further background, see \cite{lbcforbeginner}, \cite{lyubashevsky2013toolkit}, \cite{regev2010learning}. We begin with a general introduction to post-quantum cryptography (PQC).

\vspace{0.2cm}
\noindent {\bf Post-quantum cryptography.}  Cryptography relies on complex mathematical problems, such as the discrete logarithm and integer factorization problems. The RSA \cite{rsa1978rsa} and Diffie-Hellman protocols \cite{diffie2022new}, which protect our communication networks, are based on these two problems. However, rapid progress in quantum computing poses a serious threat to these foundations. The Shor algorithm \cite{shor1994algorithms}, introduced by Peter Shor in 1994, solves both problems in polynomial time on a large enough quantum computer, rendering the protocols insecure. This leads us to the field of post-quantum cryptography \cite{surveyPQC}, which focuses on cryptographic protocols that are resistant to quantum attacks. This provides an alternative for the soon-to-be-broken schemes. The main post-quantum approaches are code-based, isogeny-based, multivariate, hash-based, and lattice-based cryptography. Among these, lattice-based schemes are central for their simplicity, flexibility, and strong security guarantees. 

A problem is considered hard in the worst case if it is hard for at least one instance, and average-case hard if it is hard for most instances from a given distribution. Worst-case hardness provides theoretical assurance but may not guarantee security, as random instances might too often correspond to weak instances. To guarantee strong security, we require random (average) instances that are likely to be as hard as the worst case instances. Proving such a property is done by a process called \emph{worst-case-to-average-case reduction}.  
We define several widely studied worst-case hard lattice problems in the following section.
\bigskip

\subsection{Hard lattice problems}

We define some of the fundamental lattice problems believed to be hard in the worst case. 

\vspace{0.2cm}
\noindent \emph{Shortest vector problem (SVP):} 
Given a basis ${B}$ of a lattice ${L},$ find a shortest non-zero vector of the lattice. Explicitly, find a nonzero vector $\boldsymbol{x} \in {L}$ such that $\|\bx\|=\lambda_1({L}).$
The approximate version, \emph{approximate SVP} problem (SVP$_\gamma$), asks for a nonzero  $\boldsymbol{x}\in {L}$ such that $\lambda_1({L}) \leq \gamma(n)\|\boldsymbol{x}\|$ where $\gamma(n)\geq 1.$ The \emph{GapSVP$_\gamma$} problem (decision SVP$_\gamma$) asks to decide if $\lambda_1({L})\leq r$ or $\lambda_1({L})\geq \gamma r$ where $r \in \mathbb{Q}.$ 

\vspace{0.2cm}
 
\noindent \emph{Closest vector problem (CVP):}
Given a basis  ${B},$ of a lattice ${L},$  and a target vector $ \boldsymbol{t} \notin {L},$ find a vector in  ${L}$ that is closest to $\boldsymbol{t}.$
When $dist (t,L)\leq d$ where $d\in \mathbb{Z}^+$ we refer to this as a \emph{bounded distance decoding problem (BDD)}. The approximate version, CVP$_\gamma$, asks to find a lattice vector at distance at most $\gamma.$ And the decision version, GapCVP$_\gamma$ asks to decide whether $dist (t,\mathcal{L}) <1$ or $dist (t,{L}) <\gamma.$ 

  \vspace{0.2cm}
 \noindent\emph{Shortest independent vector problem (SIVP):} Given a lattice $L,$ find $n$ linearly independent vectors $\bv_1,\ldots,\boldsymbol{v}_n$ in $L$ such that $\max\limits_i \|\boldsymbol{v}_i\|\leq \lambda_n(L).$ The approximate version, SIVP$_\gamma$ finds these vectors with length at most $\gamma \lambda_n(L).$ On the other hand, the decision version, GapSIVP$_\gamma$ asks to determine if $\lambda_n(L(B)) \leq d$ or $\lambda_n(L(B))>\gamma d.$ 

 \vspace{0.2cm}
 
 \noindent\emph{Generalized shortest independent vector problem (GIVP$_{\gamma}^{\phi}$):} Given a lattice $L(B)$ of dimension $n,$ find $n$ linearly independent vectors $\bv_1,\ldots,\boldsymbol{v}_n$ in $L$ such that $\max\limits_i \|\boldsymbol{v}_i\|\leq \gamma\phi(L(B))$ where $\phi$ is an arbitrary real valued function of a lattice and $\gamma(n)\geq 1.$ When $\phi =\lambda_n,$ we get the SIVP$_\gamma.$ 
\vspace{0.2cm}

The hardness of these lattice problems has been widely studied due to their importance in applications. In particular, the SVP and the CVP, along with their approximate versions, form the basis for many secure post-quantum cryptographic schemes. The first results showing that CVP is computationally hard date back to Van Emde Boas \cite{van1981another}, who showed that CVP is NP-hard via a deterministic reduction. 
 That is, with high probability, any problem in NP can be reduced in polynomial time to an instance of CVP. Based on the similarities between SVP and CVP, he further conjectured SVP was also NP-hard. After almost two decades, Ajtai in \cite{ajtai1998shortest} showed that SVP with the $l_2$-norm is NP-hard for randomized reduction, thus proving the van Emde Boas conjecture. 

In practice, one typically works with the approximate variants of SVP and CVP, which are also NP-hard for a small enough approximation factor $\gamma$. Arora et al. \cite{arora1997hardness} showed that the approximate CVP is NP-hard within any constant. Micciancio \cite{micciancio2001shortest} later showed that for any $l_p$-norm, approximate SVP within any constant factor  $\gamma <\sqrt[p]{2}$ is hard under some random reductions. Specifically, in the Euclidean norm, the approximate SVP is NP-hard with any factor $\gamma <\sqrt{2}.$ Extending on \cite{arora1997hardness}, Dinur \emph{et al.} in \cite{dinur1998approximating} later showed that the approximate CVP problem in an $n$-dimensional lattice is NP-hard for a factor $\gamma=n^{\frac{c}{\log \log n}}$ for some constant $c>0.$ 

Further studies on the hardness of these problems have been conducted within an exponential factor. However, the results show that these problems are unlikely to be NP-hard. In particular, the approximate CVP and SVP are highly likely not NP-hard within a factor $\sqrt{n}$ \cite{goldreich-limits, aharonovNP}.  
Despite these limitations, these problems are computationally hard in the worst case and form the basis of modern cryptography. The connection between these worst-case hard lattice problems and security brings us to the field of lattice-based cryptography.

The study of lattice-based cryptography began with Ajtai's ground breaking work in \cite{ajtai1996}, which introduced the first average-case hard problem, \emph{the short integer solution} (SIS). On a high level, the SIS problem asks one to recover a short nonzero vector $\bz\in \mathbb{Z}^m$ given a random matrix $A \in \mathbb{Z}^{n \times m}$ satisfying $A\bz=0 \mod q.$ The requirement for $ \boldsymbol{z} \neq \boldsymbol{0}$ and short is what makes this problem difficult. This hardness was shown by a worst-case-to-average-case reduction from a well known worst-case hard problem. Building on this, in 2005 Regev introduced the learning with errors (LWE) problem \cite{Regev}, which is a noisy analog of SIS. 
 This problem plays a central role in LBC. On a high level, the LWE problem is the following: given arbitrary independent samples of noisy linear equations $(\ba,b=\ba\cdot \bs+e) \in \mathbb{Z}_q^n\times \mathbb{Z}_q$, recover the secret $\bs\in \mathbb{Z}_q^n.$ The error $e,$ is sampled from a Gaussian distribution.
The choice of the error distribution plays a critical role in the hardness of the problem as well as in the  correctness of the public key encryption scheme based on the LWE problem. We describe this below. 

\subsection{Gaussians, variational distance, and the smoothing parameter}\label{DG} 
A continuous Gaussian function centered at $c\in \mathbb{R}^n$ defined over $\mathbb{R}^n$ is given by $$ \rho_{s,c}(\boldsymbol{x})=\exp{\frac{-\pi \|\boldsymbol{x}-\boldsymbol{c}\|^2}{s^2}},$$ where $s=\sqrt{2\pi}\sigma$ with $\sigma$ being the standard deviation.
Normalizing this function by $1/s^n$ defines the corresponding probability density function
 $$g_{s,c}(\boldsymbol{x})=\frac{1}{s^n}\exp{\frac{-\pi \|\boldsymbol{x}-\boldsymbol{c}\|^2}{s^2}}.$$
  If $c=0$, we denote $\rho_{s,0}=\rho_s, g_{s,0}=g_s$. 

Let us consider the $L$-periodic function $g_{s,L}(\bx)=\sum_{\lambda \in L}g_{s,\lambda}(\bx)$, which is a probability density function when restricted to $\mathbb{R}^n/L$. 
We can now define the 
\emph{discrete Gaussian distribution}
$D_{L,s,c}$ as $$D_{{L},s,c}(\bx) = \frac{g_{s,c}(\bx)}{g_{s,L}({c})}$$ for $\bx \in {L}$ and zero otherwise. Again, we denote $D_{L,s,0}=D_{L,s}.$

\vspace{0.2cm} An important tool used in lattice-based schemes is the \emph{smoothing parameter}, which ensures that the noise parameter $s$ is large enough to hide the underlying secret structure. More rigorously, it determines the minimal noise level that guarantees indistinguishability from a uniform distribution by a maximum gap $\delta$. By ``gap'', we mean the \emph{variational distance} (equivalently, \emph{statistical distance}) between two distributions $p_X$ and $p_Y$; 
\begin{equation}\label{eq:vardist}
V(p_X,p_Y)=\int_{\mathbb{R}^n}|p_X(\boldsymbol{x})-p_Y(\boldsymbol{x})|d\boldsymbol{x}.
\end{equation}

\vspace{0.1cm}
\noindent\emph{Smoothing parameter \cite{micciancio2007worst}:} Let $L \subset \mathbb{R}^n$ be a full-rank lattice and $\delta > 0$.
The smoothing parameter $\eta_\delta(L)$ is defined as
\begin{equation}\label{eq:smooth}
\eta_\delta(L)= \inf \left\{ s> 0 \;\middle|\;\rho_{1/s}(L^* \setminus \{0\}) \;\le\; \delta \right\}.
\end{equation}

In other words, fixing $s$ to be the infinimum ensures that the variational distance between the lattice Gaussian and the uniform distribution on the Voronoi cell is at most $\delta$.

The smoothing parameter is closely related to the flatness factor, which we will define in Section \ref{sec:IT-sec}, Eq. \eqref{eq: flatness factor definition}. 

 \begin{rem}
   The continuous Gaussian is used for analysis such as defining the smoothing parameter and for reduction proofs. In contrast, the discrete Gaussian is used for the actual cryptographic constructions, such as sampling the errors.  We will use both in this work as needed depending on the context.  
 \end{rem}
\subsection{Learning with errors (LWE) and its variants}
We define the LWE problem and some of its variants like the ring learning with errors (RLWE), the polynomial learning with errors (PLWE) problems and the module learning with errors problem (MLWE).  We closely follow the definitions in \cite{Regev}, \cite{LPR2010}, \cite{BV}.
\leavevmode\\[0.2cm]
{\bf Learning with errors (LWE) problem.}
Let $n\geq 1$ and $q=q(n) \geq 2$ be integers. Let $\bs \in \mathbb{Z}_q^n$ be a secret vector sampled uniformly where $n$ is the security parameter. Additionally, let $\chi$ be the error distribution that follows a discrete Gaussian sampled over $\mathbb{Z}$ and reduced modulo $q.$
\leavevmode\\[0.1cm]
 \noindent \emph{LWE distribution:}
The LWE distribution $\mathcal{A}_{\bs,\chi}$ over $\mathbb{Z}_q^n \times \mathbb{Z}_q$ is defined as follows: sample 
$ \ba\leftarrow \mathbb{Z}_q^n$ uniformly, $e \leftarrow \chi,$ and output $(\ba,b) \in \mathbb{Z}_q^n \times \mathbb{Z}_q$ where $ b=\langle\bs,\ba\rangle + e \pmod q.$ The addition is performed modulo $q.$
\leavevmode\\[0.1cm]
 \noindent\emph{Search LWE (LWE$_{q,\chi}$):}
Given an arbitrary number of independent samples $(\ba_i,b_i),$ drawn from the LWE distribution $\mathcal{A}_{\bs,\chi},$ the search LWE problem asks to find the secret vector $\bs.$
\leavevmode\\[0.1cm]
 \noindent \emph{Decision LWE:}
Given arbitrary many independent samples the decision problem asks to determine with non-negligible advantage whether the sample $(\ba_i,b_i),$ is from the LWE distribution $\mathcal{A}_{\bs,\chi}$ or from a uniform distribution.

In matrix form:  $(A, \bb) \in \mathbb{Z}_q^{n \times m} \times \mathbb{Z}_q^m$ where $A$ is a matrix whose columns are the vectors $\ba_i\in \mathbb{Z}_q^n$ and $\bb^t =\bs^tA+ \be^t \pmod q$ with $m$ the number of samples. Let 
\[
L=\{\boldsymbol{y}\in \mathbb{Z}^m:A^T\bz=\boldsymbol{y} \pmod q \;\text{ for some }\bz\in \mathbb{Z}^n\}.
\]
Observe that $L$ is a full rank integer lattice in $\mathbb{R}^n$, when we choose $m=n$ linearly independent $\ba_i$'s. The LWE problem can be viewed as a bounded distance decoding problem on $L$ where $\bb$ is a closest vector to a lattice point.
 
The relevance of LWE lies in its reduction from a worst-case hard problem. We state the hardness result below.

\begin{thm}(\cite{Regev}, Theorem 1.1)
Let \( n,q \) be integers, $\alpha \in [0,1)$ be such that $\alpha q>2\sqrt{n}$ and $\chi$ a Gaussian distribution. If there exists an efficient algorithm that solves LWE$_{q,\chi},$ then there exists an efficient quantum algorithm that approximates the decision version of SVP (GapSVP) and SIVP to within $\tilde{{O}}(n/\alpha)$\footnote{The function $f(n) = {O}(g(n))$ if there exist constants $c>0$ and $n_0$ such that 
$f(n) \le c\,g(n)$ for all $n \ge n_0$. We denote $\tilde{O}(f(n))=O(f(n)poly\log(n)).$} in the worst case. 
\end{thm}

The LWE decision and search versions are equivalent when the integer modulus $q$ is prime \cite{Regev}. 
Although LWE offers strong security guarantees, LWE-based schemes have a quadratic overhead in the key size in terms of the security parameter $n$, rendering them inefficient for practical purposes. To address this concern, Lyubashevsky, Peikert, and Regev \cite{LPR2010}, \cite{lyubashevsky2013toolkit} introduced the ring variant of LWE, known as the ring learning with errors (RLWE) problem, which only induces a linear overhead. This variant was originally formulated in the so-called dual form. However, the primal version is preferable in practice. Since the dual and primal formulations are known to be equivalent \cite{RSW2018}, we restrict ourselves to the primal version, which we refer to as RLWE. We closely follow the definition in \cite{EHL2014}.  

\vspace{0.1cm}
\noindent {\bf Ring learning with errors (RLWE) problem.} Let $n\geq 1$ and $q=q(n)\geq 2$ be an integer modulus. Let $K$ be a number field of degree $n$ and $R=\mathcal{O}_K$ its ring of integers. 
Set $R_q=\mathcal{O}_K/q\mathcal{O}_K$ 
and let $\chi$ be the discrete Gaussian distribution obtained by sampling over the ideal lattice $L=\Sigma(R)$ (see Section~\ref{constr}). Let $s \in R_q$ be sampled uniformly at random.
\vspace{0.1cm}

\leavevmode\\[0.1cm]
\noindent\emph{RLWE distribution ($\mathcal{A}_{{s},\chi}$):}
The RLWE distribution $\mathcal{A}_{s,\chi}$ in $R_q \times R_q$ is defined by uniformly sampling
$ a \leftarrow R_q,$
   $e \leftarrow \chi,$ and returning
    $(a,b) \in R_q \times R_q$ where $b=a s + e \mod q.$ 
\leavevmode\\[0.1cm]
\noindent \emph{Search RLWE (RLWE$_{q,\chi}$):}
For an arbitrary number of independent samples $(a_i,b_i),$ drawn from the RLWE distribution $\mathcal{A}_{s,\chi},$ the RLWE$_{q,\chi}$ problem asks to recover the secret $s.$ 
\leavevmode\\[0.1cm]
\noindent \emph{Decision RLWE (D-RLWE$_{q,\chi}$):}
Given an arbitrary number of independent samples $(a_i,b_i),$ the D-RLWE$_{q,\chi}$ asks to determine with non-negligible advantage whether it came from the RLWE distribution $\mathcal{A}_{{s},\chi}$ or a uniform distribution. 

The following result guarantees the hardness of this problem.
\begin{thm}[\cite{LPR2010}] Let $K=\mathbb{Q}(\zeta_m)$ be the $m$th cyclotomic number field with degree $n=\varphi(m)$ and $R=\mathcal{O}_K$ be its ring of integers. Let $\alpha=\alpha
(n)$ and let $q=q(n) \geq 2,$ $q\equiv 1 \pmod m$ be a poly(n)-bounded prime such that $\alpha q>\omega(\sqrt{\log n})$\footnote{The function $\omega(f(n))$ grows asymptotically faster than $f(n)$.}. Then there is a polynomial time reduction from $\tilde{O}(\sqrt{n}/\alpha)$-approximate SIVP (or SVP) on ideal lattices in $K$ to the problem of solving D-RLWE$_{q,\chi}$
\end{thm}

The assumption that $K$ is a cyclotomic field is only needed in the reduction from RLWE$_{q,\chi}$ to D-RLWE$_{q,\chi}$. 

This reduction has been extended to any Galois extension \cite{EHL2014}. 
Unlike LWE-based schemes, cryptographic protocols based on RLWE have linear overhead in the degree of the number field due to the added structure \cite{lyubashevsky2013toolkit}. However, this added structure may introduce weaknesses that are not present in plain LWE \cite{ELOS}.

One way to improve the security is by increasing the field extension degree, thereby also increasing the degree of the involved polynomials. This has an obvious adverse effect on the computational complexity. A better alternative is by another structured LWE variant, the general or module LWE framework of \cite{BGV}, later formalized as MLWE with a worst-case hardness reduction by Langlois and Stehl\'e \cite{langlois2015worst}.  We will closely follow \cite{boudgoust2023hardness,langlois2015worst}. Essentially, in a sample $(a,b)$ we now  consider $a$ to be a vector of length $d$ instead of a single ring element or polynomial ($d=1$), resulting in a module structure over a ring. This enables increasing the security level by increasing the module rank $d$, while keeping the underlying field extension degree and ring intact. Moreover, instead of a single distribution, we will now need a family of distributions to match the varying $d$.

\vspace{0.2cm}
\noindent {\bf Module learning with errors (MLWE) problem.} Let $K$ be a number field of degree $n,$ $R=\mathcal{O}_K,$  $\Psi$ a family of distributions on $K_{\mathbb{R}}$ and the torus  $\mathbb{T}=K_{\mathbb{R}}/R$. For $q,d$ positive integers with $q\geq 2$ and $d\geq 1$,
let $\bs\in R_q^d$ be the secret and $\psi \in \Psi.$ Let $N = nd$ denote the dimension of the corresponding module lattice. We define the primal version of the problem as presented in \cite{boudgoust2023hardness}.

\vspace{0.2cm}

\noindent\emph{Module learning with errors distribution.} The MLWE distribution $\mathcal{A}_{s,\psi}^{\mathcal{M}}$ is obtained by sampling $\ba \leftarrow \mathcal{U}(R_q^d),$ $e\leftarrow \psi$ and returning $(\ba,b) \in R_q^d \times \mathbb{T}$ where $b=q^{-1}\langle \ba, \bs\rangle + e \pmod R.$

\vspace{0.2cm}

\noindent\emph{Search MLWE:} Given arbitrary many samples from $\mathcal{A}_{\bs,\psi}^{\mathcal{M}},$ the search MLWE problem, MLWE$_{q,\Psi}$, asks to recover $\bs.$

\vspace{0.2cm}

\noindent\emph{Decision MLWE:} Let $\Upsilon$ be a distribution on a family of distributions on $K_{\mathbb{R}}.$ The decision MLWE, D-MLWE$_{n,d,q,\Upsilon},$ is to distinguish with non-negligible advantage between arbitrary many independent samples from $\mathcal{A}_{\bs,\psi}^{\mathcal{M}}$ and the same number of independent samples from $\mathcal{U}(R_q^d\times \mathbb{T}).$

The MLWE problem was originally motivated by the goal of building a fully homomorphic scheme without bootstrapping. Langlois \emph{et al.} in \cite{langlois2015worst} later proved that MLWE  is as hard as solving approximate SIVP over module lattices in the worst case. We refer to their paper for a thorough discussion on the security of MLWE.

\begin{thm}
  Let $\varepsilon(N) = N^{-\omega(1)}$, $\alpha \in (0,1)$, and $q \geq 2$ of known factorization such that
\[
\alpha q > 2\sqrt{d} \cdot \omega(\sqrt{\log n}).
\]
There is a quantum reduction from solving GIVP$_{\gamma,\varepsilon}^{\eta}$ over module lattice in polynomial time (in the worst case, with high probability) to solving MLWE$_{q,\Psi_{\leq \alpha}}$ in polynomial time with non-negligible advantage, where
\[
\gamma = \frac{8N d \cdot \omega(\sqrt{\log n})}{\alpha}.
\]

Assume that $q$ is prime, $q \leq \mathrm{poly}(N)$, and that $q \equiv 1 \pmod{m}$ where $m$ is the conductor of a cyclotomic field. Then there exists a polynomial-time reduction from MLWE$_{q,\Psi_{\leq \alpha}}$ to D-MLWE$_{q,\Upsilon_{\alpha}}$.   
\end{thm}

 As mentioned earlier, the module learning with errors is a generalization of the previous variants of LWE. More precisely, setting $n=d=1$ corresponds to the basic LWE, while $d=1, n>1$ yields RLWE.

Although this is more efficient than standard LWE, in practice, concrete polynomial rings are often preferred for their practical advantages and simplicity. This leads us to the polynomial learning with errors (PLWE) problem, described in the following section.
\leavevmode\\[0.2cm] 
\noindent {\bf Polynomial learning with errors (PLWE)  problem.} Let $n\geq 1$ and $q=q(n)\geq 2.$ Set $R(f)=\mathbb{Z}[x]/(f(x))$ to be the polynomial ring  and $R_q(f)=R(f)/qR(f)$ and $\chi$ be a  discrete Gaussian over $\rho_f(R(f))$ (see Section~\ref{constr}). Let $s \in R_q$ be sampled uniformly at random.
\leavevmode\\[0.1cm]
\noindent\emph{PLWE distribution ($\mathcal{B}_{f,s,\chi}$):}
The PLWE distribution $\mathcal{B}_{f,s,\chi}$ in $R_q(f) \times R_q(f)$ is obtained by uniformly sampling
$ a \leftarrow R_q(f),$
   $e \leftarrow \chi,$ and returning
    $(a,b) \in R_q(f) \times R_q(f)$ where $ b=as + e \pmod q.$ 
\leavevmode\\[0.1cm]
\noindent \emph{PLWE (Search PLWE$_{f,q,\chi}$):}
For an arbitrary number of independent samples $(a_i,b_i),$ drawn from the PLWE distribution $\mathcal{B}_{f,s,\chi},$ the PLWE$_{f,q,\chi}$ search problem asks to find the secret $s.$ 
\leavevmode\\[0.1cm]
\noindent \emph{Decision PLWE (D-PLWE$_{f,q,\chi}$):}
Given an arbitrary number of independent samples $(a_i,b_i),$ the D-PLWE$_{f,q,\chi}$ asks to determine whether it was drawn from the PLWE distribution $\mathcal{B}_{f,s,\chi}$ or a uniform distribution.  

Unlike the RLWE problem, which admits an worst-case-to-average-case reduction, the PLWE problem has such reductions only for powers of two cyclotomic fields. A natural step to extend this reduction to a broader class of fields will be to study the equivalence between RLWE and PLWE. This allows us to enjoy both efficiency and security advantages. In the following section, we define the equivalence between RLWE and PLWE and survey known results on this topic.

\subsection{Equivalence between RLWE and PLWE}  

 The RLWE and PLWE problems are said to be \emph{equivalent} if there exists an algorithm that transforms a RLWE sample into a PLWE sample and vice versa in polynomial time, incurring a noise increase that is polynomial in the degree of the number field. 
This sample transformation is performed using the following map. 
$$V_f: \mathbb{Z}[x]/(f(x)) \to \sigma_1(\mathcal{O}_K)\times\dots\times\sigma_n(\mathcal{O}_K)$$
\begin{equation}\label{trans}
    \sum_{i=0}^{n-1} a_i x^i 
    \mapsto \underbrace{\begin{bmatrix} 1 & \theta_1 & \dots & \theta_1^{n-1} \\ 1 & \theta_2 & \dots & \theta_2^{n-1} \\ \vdots & \vdots & \ddots & \vdots \\ 1 & \theta_n & \dots & \theta_n^{n-1} \end{bmatrix}}_{V_f} \begin{bmatrix} a_0 \\ a_1 \\ \vdots \\ a_{n-1} \end{bmatrix}
\end{equation}

This transformation incurs some noise distortion, which is quantified by the condition number of $V_f$ \cite{RSW2018},  defined by $\mathrm{Cond}(V_f)=\|V_f\|\|V_f^{-1}\|,$ where $$\|V_f\|= \sqrt{\mathrm{Tr}(V_f V_f^*)}.$$ This reduces the question of equivalence to analyzing the condition number of the transformation matrix $V_f.$ This topic has been widely studied and we highlight some known results.

Ducas and Durmus in \cite{ducas2012ring} showed that equivalence holds for power-of-two cyclotomic fields where the change-of-basis matrix is a scaled isometry. In \cite{RSW2018}, equivalence was established for a certain ad hoc class of fields. Further results were obtained for cyclotomic fields whose conductors are divisible by two distinct prime factors \cite{blanco2022rlwe}, \cite{scala2020condition}. This result was later extended in \cite{BARB2025} to cyclotomic fields whose conductors are divisible by six distinct prime factors. Later in \cite{NotEqui} it was shown that equivalence fails for general cyclotomic fields. 

Other classes of fields have also been studied, such as the maximal real subfields of cyclotomic fields \cite{ahola2025fast, bolanos2025fast} and the cyclo-multiquadratic fields (the composition of cyclotomic and multiquadratic fields) \cite{BARB2025}. 

\vspace{0.1cm}

The security of the above structured lattice problems is based on our choice of $f$, the modulus polynomial, and the defining polynomial of the number field. A wrong choice of $f$ could make the scheme vulnerable to attacks. We review some of the attacks that exploit these additional structures. For a more thorough survey we refer to \cite{nchiwo2026cryptanalysis}.

\subsection{Cryptanalysis of RLWE/PLWE} In this section, $R_q=\mathbb{F}_q[x]/(f(x))$ where $f(x)$ is a monic irreducible polynomial in $\mathbb{Z}[x]$ and $q$ prime. 

An attack on PLWE that exploits the algebraic structure of $R_q,$ was first introduced by Eisenträger, Hallgren, and Lauter in \cite{EHL2014}. Let $\alpha$  be a root of $f(x)$ i.e $f(\alpha)=0 \mod q.$ They showed that when $\alpha=1$, an attacker, given arbitrary PLWE samples $(a_i, b_i)\in R_q^2$, can efficiently distinguish them from uniform samples. In the same work, they extended this idea to the case where $\alpha$ has a \emph{small multiplicative order} $r \mod q.$ Later, Elias, Lauter, Özman, and Stange \cite{ELOS}, further constructed a  distinguishing 
attack on PLWE when $\alpha$ has a \emph{small residue}.

Building on these root-based attacks, the authors of \cite{blanco2023trace} constructed a similar attack that makes use of the number-theoretic properties of the trace function without requiring the order of the root to be small. Specifically, they show that if $f(x)$ has a quadratic factor whose root has trace zero, then the adversary can identify whether the given samples are uniform or PLWE. This idea was subsequently generalized in \cite{blanco2025generalized}, which extended the result from quadratic factors to factors of higher degree: it showed that if $f(x)$ contains a factor of the form $x^n+\rho$, where $\rho$ is an element such that the root has trace zero, then a similar distinguishing attack applies. Furthermore, using a similar strategy of exploiting roots with trace zero, another attack was designed in \cite{barbero2023cryptanalysis} 
on PLWE over a subring $R_{q,0}\times R_q$ of a cyclotomic field, where the integer modulus $q$ does not split completely.

Another attack proposed in \cite{EHL2014} and addressed in \cite{smearingattack} is the smearing attack. This distinguishing attack is performed by analyzing the behavior of the error to distinguish them from uniform samples.

The natural question that arises is: can some of these attacks on PLWE be extended to RLWE? This question has been answered in the affirmative for some. A natural step in this direction will be to transform RLWE samples into PLWE (see Equation \eqref{trans} ) and then apply one of the listed attacks to the PLWE samples. 
In \cite{ELOS}, the authors demonstrate how the attack described in \cite{EHL2014} can be applied to the RLWE decision problem, given that some conditions are satisfied. This result was later extended in \cite{castryck2016provably} to the RLWE search problem, achieving $100\%$ success probability with fewer samples.

It is worth mentioning that these known attacks do not threaten the
NIST-standardized schemes. However, these attacks confirm why we should stick to the parameters of the standardized schemes. They also provide a list of parameters to watch for when constructing new schemes or improving standardized ones. 

\subsection{Further applications of lattice-based cryptography}
One interesting property of lattice-based cryptographic (LBC) schemes is their flexibility, which enables them to be applied across a wide range of applications. We briefly describe two notable examples below: homomorphic encryption and private information retrieval (PIR). We conclude this section by highlighting the new standardized schemes based on lattices.

\leavevmode\\[0.1cm]
\emph{\bf Homomorphic encryption.} This is a form of encryption that allows computations to be done on encrypted data without decrypting it. The decrypted result matches the result of the same operation performed on the original plain data. There are three main types: partially homomorphic encryption (PHE), which allows for a single type of operation (either addition or multiplication) on encrypted data; somewhat homomorphic encryption (SHE), which allows for a limited, but greater than one, number of operations on encrypted data; and fully homomorphic encryption (FHE), which allows for an unlimited number of either of the two operations. Homomorphic encryption plays an important role in, e.g., cloud computing, secure voting, and medical data analysis.

The problem of constructing a fully homomorphic encryption scheme was first introduced by Rivest, Adleman, and Dertouzos in 1978 \cite{rivest1978data}. For more than three decades, the existence of a solution remained an open problem. During this period, partial homomorphic schemes like RSA \cite{rsa1978rsa} and  ElGamal \cite{elgamal1985public} enabled unlimited modular multiplications, while others allowed for unlimited modular arithmetic \cite{paillier1999public}, \cite{benaloh1994dense}. 

Gentry, in \cite{gentry2009thesis}, \cite{gentry2009fully}  gave the first construction of an FHE scheme based on lattice-based cryptography. This FHE scheme supports both addition and multiplication operations on ciphertext, for an arbitrary number of computations. To achieve this, he starts from a somewhat homomorphic encryption scheme, which has limitations due to the noise growth when we add or multiply encrypted data. This noise growth may lead to decryption errors if the noise becomes too high. Addressing this, he shows that in any bootstrapped scheme, this SHE scheme can be converted into an FHE scheme. The security of this scheme is partly based on worst-case problems over ideal lattices. Furthermore, in \cite{BGV}, this result is extended by removing the ideal lattice condition. Over the years, other schemes based on lattice-based problems have been proposed that offer improved efficiency. 

One of the well known applications of homomorphic encryption is the single-database computationally-Private Information
Retrieval (cPIR), which we introduce below.
 
\vspace{0.3cm}

\noindent\emph{\bf Private information retrieval (PIR)}  Introduced by Chor \emph{et al.} in \cite{Chor}, PIR allows a user to retrieve an item from a storage system or cloud in possession of a database without revealing which item is retrieved. The main idea is as follows: given a system holding a database consisting of a set of 
 elements $D_1,\dots,D_m,$ retrieve the $i$th element $D_i$
 without revealing $i$ to the system owner. A
 naive solution will be to retrieve the entire database and  discard all entries but  the one of interest. However, this will be at a cost proportional to the number of data items $O(m)$, and therefore not practical when the database is large. It has been shown that this is the only way to guarantee information-theoretic privacy if we store the database on one server. However, using \emph{multi-server schemes} we can do much better. 

Here, the user sends masked queries to different non-colluding servers, ensuring that no subset of servers of size below a design threshold learns the original query. The methods utilized in this approach typically draw from coding theory, and various extensions to the problem have been made. The literature is vast and as our core topic is lattices, we refrain from expanding our reference list by these coding-theoretic works and simply refer to the brief tutorial \cite{bits} and the references therein. 

Another way is to only use a single server but instead rely on computational security and hence cryptographic techniques to hide the query from the server, originally introduced in \cite{single-server} where a scheme based on quadratic residues was constructed. This method is generally known as \emph{single-server computationally private information retrieval (cPIR).} The drawback of the original cPIR scheme is that it is computationally more expensive than the naive method of downloading everything and therefore not practical \cite{single-server}. This aspect has later been improved by several works, notably in \cite{XPIR} where, using LBC, much more efficient cPIR schemes based on LWE and RLWE were constructed. 

\vspace{0.3cm}

\noindent\emph{\bf NIST standardized schemes.} Due to the fast progress in the field of quantum computing, the National Institute of Standards and Technology (NIST) launched a PQC standardization process (competition) in $2016$ aimed at identifying quantum-safe cryptographic protocols. A call for proposals \cite{NistAnnouncement} was made in which researchers were invited to submit candidates which were evaluated at several rounds based criteria that included security, efficiency, and practical implementation. A total of $82$ algorithms were submitted, of which $69$ were accepted into the first round\cite{NistRound1}, $26$ advanced to the second round \cite{NistRound2}, 
and $15$ selected as third round finalist \cite{NistRound3}. Of these candidates, $4$ was chosen for standardization \cite{NistDecision} and the rest moved to the fourth round \cite{NistRound4} and are currently undergoing further analysis. Among the four schemes chosen for standardization, CRYSTALS-Kyber (ML-KEM), CRYSTALS-Dilithium (ML-DSA), and Falcon are based on hard lattice-problems. In particular, Kyber and Dilithium are based on hard problems over module lattices \cite{avanzi2019crystals} such as MLWE. 
It is worth noting that the above listed attacks on RLWE/PLWE do not threaten any of these standardized schemes.

\section{Lattice codes for wireless security}\label{sec:IT-sec}

In the previous section, we have seen how lattices can be applied to provide computational security due to several hard lattice problems, such as the shortest vector problem. In addition to computational security where we assume that the adversary has bounded computational resources, lattices can also be used to provide \emph{physical layer security} (PLS), which for its part relies on the concept of \emph{information-theoretic security}.  Here, the adversary can have unlimited computational power, since the security is based on (full or partial) lack of relevant information. 

While traditional cryptographic methods like AES (the Advanced Encryption Standard) or RSA operate at higher protocol layers, PLS provides a complementary layer of defense. As the 6th generation (6G) communication networks move toward high-mobility, low-latency, and decentralized networks (e.g., Internet of Vehicles), traditional key exchange mechanisms can become bottlenecks. Lattice-based PLS offers a way to establish security guarantees by leveraging the intrinsic physical properties of the wireless medium, making it a vital research area for future secure communication systems against both classical and potential quantum threats. We refer to the recent white paper \cite{whitepaper} for a general introduction to the utility of physical layer security. 

Let us start by defining relevant notions in information theory. 

\subsection{Basic notions in information theory and related security paradigms}

For a general reference for information theory and information-theoretic security, we refer to \cite{cover_IT, oggier-tutorial}. Let $X$ and $Y$ be two discrete random variables taking values from respective sets $\mathcal{X}$ and $\mathcal{Y}$. The \emph{entropy} of $X$ is
$$
H(X)=-\sum_{x\in\mathcal{X}}p(x)\log p(x),
$$
where $p(x)$ is the probability of $x$.  The entropy of a continuous random variable is defined analogously. 

Let us denote the conditional entropy of $X$ given $Y$ by $H(X|Y)$. The \emph{mutual information} can then be defined as 
$$I(X;Y)=H(X)-H(X|Y).$$

Assume now that we are sending a message $X$ and the adversary is gaining access to $Y$. The secrecy can be measured by how much information can be obtained based on $Y$, quantized by \emph{information leakage} $I(X;Y)$. The mutual information is zero, yielding \emph{perfect secrecy}, if and only if $X$ and $Y$ are independent. In cryptography, typically operating over positive characteristic and finite structures, this can be achieved by adding uniformly random noise to $X$, referred to as one-time pad. On a wireless (physical) channel, however, the noise is not \emph{selected} by the user as in cryptography but is \emph{induced} by the physical conditions of the channel and the equipment used and as such largely beyond our control. Such noise is typically real or complex Gaussian, and we cannot (non-asymptotically) obtain perfect secrecy anymore. Instead, the goal in physical layer security is to minimize information leakage while guaranteeing good decoding probability for the legitimate user. 
This will be our focus in the rest of this section. We refer the interested reader to \cite{viterbo,belfiore,chorti,costa2017lattices} and references therein for more details on lattice-based reliable and secure wireless communications. 

\subsection{Wiretap channels and lattice coset codes} 

Lattice codes, defined as finite collections of lattice vectors within a bounding region, are effective tools for wireless communications. The channel model for single-input single-output (SISO) is described as
$$
\mathbf{y}=H\mathbf{x}+\mathbf{n}\in \mathbb{R}^n,
$$
where $\mathbf{x}\in L$ is the message, $\mathbf{n}\in \mathbb{R}^n$ is additive white Gaussian noise (AWGN), and $H\in \mathbb{R}^{n\times n}$ represents random fading. For an \emph{AWGN channel}, $H=I_n$, while for a \emph{Rayleigh fast fading channel}, $H$ is a diagonal matrix with independent Rayleigh distributed entries. Maximum-likelihood (ML) decoding here is equivalent to the closest vector problem in a (distorted) lattice. This may sound counter-intuitive after the previous section on lattice-based cryptography, where it was pointed out that many lattice problems, including the CVP, are computationally hard. Fortunately for us, now the security is based on information-theoretic notions, not on computational hardness. This means that, in practice, we can resort to relatively low-dimensional lattices making the CVP feasible for the legitimate receiver. However, in order to approach the theoretical perfect secrecy capacity \cite{oggier-seccap}, high-dimensional lattices are still needed.

In a channel with not too much fading and noise with respect to the signal power, measured by \emph{signal-to-noise ratio (SNR)}, reliability is optimized by maximizing the \emph{modulation diversity} $\ell:=\min_{0\neq\mathbf{x}\in L}|\left\{i: x_i\neq 0\right\}|$ and the \emph{minimum product distance} $$d_{p,min}(L):=\inf_{0\neq\mathbf{x}\in L}\prod_{i=1}^{n}|x_i|$$ for a full diversity lattice ($\ell=n$) \cite{viterbo}. At low SNR, the minimum distance $\lambda_1(L)$ dominates the decoding performance. Note that the minimum product distance coincides with the multiplicative norm defined in \ref{sec:lattice}.

Wyner's coset coding, introduced by Aaron Wyner in 1975 \cite{wyner, Wyner-Ozarow} for the \emph{wiretap channel}\footnote{Here, in addition to the legitimate receiver, we have an eavesdropper who is ``tapping the wire'' in order to intercept secret messages. With slight abuse of language, we also call a wireless channel with an eavesdropper a wiretap channel, even though there is no wire.}, is a foundational technique in information-theoretic security that achieves secure communication without relying on unproven computational assumptions or shared cryptographic keys. The original core mechanism involves partitioning a standard error-correcting code into distinct, non-overlapping subsets, \emph{i.e.}, (cosets), where each coset corresponds to a specific secret message. To transmit a message, the sender selects the appropriate coset and deliberately introduces structured randomness by transmitting a randomly chosen codeword from the coset. This approach is useful in physical-layer security as it exploits the differences in the channel quality between the legitimate user and an eavesdropper. A legitimate receiver with a stronger channel (\emph{i.e.}, lower noise) can successfully decode to identify the correct coset and recover the message, while an eavesdropper on a noisier channel is overwhelmed by the errors. Due to the random codeword selection, the eavesdropper's degraded signal lacks enough information to even determine which coset was used, mathematically ensuring that the intercepted data reveals (practically) zero information about the original message regardless of the attacker's computing power. We refer to \cite[Ch. 3]{okkothesis} and the references therein for a nice exposition on coset codes and wiretap channels for error-correcting codes, as well as their intimate connection to \emph{(homomorphic) secret sharing}. 

In our lattice code context, instead of quotients of vector spaces we consider quotients of lattices, both of which have an analogous additive group structure. A message $\mathbf{m}$ is masked by a random sublattice vector $\mathbf{r} \in L_s \subset L$, and $\mathbf{x}=\mathbf{m}+\mathbf{r}\in L/L_s$ is transmitted. Again, assuming the eavesdropper experiences relatively stronger noise, the legitimate receiver decodes reliably while the eavesdropper gains negligible information. For more details and examples, we refer to \cite{Oggier-Sole-Belfiore}. 

Next, let us make the ``negligible information'' more rigorous.

\subsection{The flatness factor and connections to theta functions}

The \emph{flatness factor} $\mathcal{E}_{L} (\sigma)$  bounds from above the deviation of the lattice Gaussian $g_\sigma$ from uniformity on the Voronoi cell $\mathcal{V}(L)$ (cf. variational distance, Eq. \eqref{eq:vardist}):
\begin{equation}
\label{eq: flatness factor definition}
\mathcal{E}_{L} (\sigma) : = \max_{\mathbf{x} \in \mathcal{V}(L)} \left| \frac{g_{\sigma,L}(\mathbf{x})}{1/\Vol(L)} - 1 \right|.
\end{equation}
The flatness factor bounds from above both the eavesdropper's correct decoding probability and information leakage \cite{Belfiore-flatness,luzzi_isit16,damir-2}, and minimizing it makes the channel appear ``flatter'' (more uniform) to the eavesdropper. 

The flatness factor can be related to the primal and dual theta series (via the Poisson summation formula) as follows:
\begin{eqnarray}
\label{eq: dual theta flatness factor formula}
  &&\Vol(L) g_{\sigma,L}(\bx) - 1 \nonumber \\
  &\leq&  \Vol(L) g_{\sigma,L}(0) - 1 \nonumber \\
&=& \frac{\Vol(L)}{(\sqrt{2 \pi} \sigma)^n} \Theta_{L} (e^{-1/2 \sigma^2}) - 1 \\
&=& \Theta_{L^*} (e^{- 2 \pi \sigma^2}) - 1 \nonumber \\
&=& \mathcal{E}_{L} (\sigma) \nonumber
\end{eqnarray}
where $\sigma^2$ is the noise variance.

A practical complication now arises. Namely, if we want to compare different lattices (of the same volume) and compare their flatness factors, how do we efficiently compute the theta function, known to be notoriously hard? In a relatively low dimension, one can use truncations and brute force point enumeration or, more efficiently, resort to the theta function approximation proposed in \cite{Amaro_approx}. However, as mentioned by the authors, the approximation is not universally good and gets worse with growing dimension or if the lattice is very skewed.

To make a connection to lattice-based cryptography, we point out that the flatness factor and the smoothing parameter (cf. Eq. \eqref{eq:smooth}) are closely related \cite{Belfiore-flatness}. To this end, let us slightly redefine the smoothing parameter up to constants by changing the variable $s$ to $\sigma=s/{\sqrt{2\pi}}$:
$$
\eta_\delta(L)=\inf\{\sigma>0\,|\, \sum_{0\neq\lambda\in\Lambda^*}e^{-2\pi^2\sigma^2||\lambda||^2}\leq \delta\}.
$$
Then we have
$$
\mathcal{E}_{L} (\eta_\delta(L))=\delta. 
$$
Essentially, if we compare two lattices (of the same volume), the one with the smaller flatness factor allows for adding smaller noise while maintaining the same variational distance, helping the correct decoding/decryption of the legitimate receiver. 

\subsection{Well-rounded lattices as theta minimizers}
Well-rounded lattices were studied and constructed for the SISO wiretap channel in \cite{Gnilke,spawc}. Interestingly, it has been shown that the minimizer of the flatness factor (equivalently, the theta function) is well-rounded \cite{damir-2}. Furthermore, it is conjectured that the closer the lattice is to being stable, the smaller the flatness factor \cite{jesse-thesis}. 

 Since they also maximize the minimum product distance \cite{damir-3} and support dense packings \cite{sarnak,coulangeon,delone}, WR lattices are excellent candidates for secure, reliable communication. As minimizing the theta function globally is difficult, \cite{camilla} focused on generic WR lattices --- increasing the first minimum and decreasing the kissing number reduces the dominating term, motivating the study of dense GWR lattices for physical layer security. 

\subsection{Related topics and generalizations}
Here, we have concentrated on the single-antenna channel model. For multi-antenna wireless communications, so-called \emph{space-time lattice codes} based on \emph{cyclic division algebras and their maximal orders} can be used \cite{OggierCDA, sethudiv, camimax, camidense}. Lattice coset codes and related design criteria for such a \emph{multiple-input multiple-output (MIMO)} wiretap channel have also been considered \cite{mirghasemi,luzzi}. The utility of well-rounded lattices for MIMO channels have been demonstrated in \cite{Gnilke-Barreal,Barreal-Karrila-Karpuk-Hollanti}. Furthermore, an analogous design criterion to minimize the lattice theta series in the $\ell^1$ (taxicab) norm instead of the euclidean norm was proposed in \cite{niklas_l1}, and related kissing number problems in \cite{niklas_kiss}.

The setting of point-to-point communications can be extended to \emph{relay channels}, where the message is relayed by an intermediate node. Lattice coset codes also come into play here via \emph{physical layer network coding}, the security of which has been considered in \cite{Navin_SPLNC} for the so-called compute-and-forward channels. Generalized theta series was considered in \cite{maiara-gentheta}, motivated by connections to the identification of stable lattices, the  \emph{lattice isomorphism problem}, and the so-called isodual \emph{secrecy gain conjecture}. The secrecy gain is closely related to the flatness factor, and measures how much better the coding lattice used is with respect to ``no coding'', \emph{i.e.}, the $\mathbb{Z}^n$ lattice, by looking at the ratio of the respective theta functions called the \emph{secrecy function} \cite{Belfiore2010Secrecy-gain, Oggier-Sole-Belfiore}. Its ultimate goal is the same as that of the flatness factor minimization: to minimize the the theta function of the eavesdropper's lattice (the reader can think of this as the ``noise'' lattice, which is used to confuse Eve).  While the flatness factor directly bounds the eavesdroppers correct decoding probability and information leakage, hence making comparison to the integer lattice obsolete, the \emph{secrecy gain} does allow the use of somewhat different type of analytical tools by looking at the maximum of the secrecy function. For more details on the secrecy gain and various conjectures related to it, we refer to \cite{Belfiore2010Secrecy-gain, Oggier-Sole-Belfiore, ame, maiara}. Finally, we mention that the maximization of the flatness factor has been studied in \cite{maiaraFF}.

For a more general and broader introduction to many of the topics discussed in Sections \ref{sec:lbc} and \ref{sec:IT-sec}, we refer the reader to the following PhD theses:  \cite{taothesis,rahithesis,amarothesis,okkothesis}.

\section{Conclusion and open problems}\label{sec:summary}

In this survey, we have discussed the theory of structured Euclidean lattices and its applications --- a very active area of current research. A great deal of research here is motivated by the original discrete optimization problems on lattices, such as packing, covering, and kissing number problems. This being said, the field of potential work here is much wider than suggested by these classical problems, both in terms of theory and applications. We mention here some additional directions for future work.
\medskip

\noindent
{\bf Minkowski conjecture} has been proved in dimensions $n \leq 10$, and the related Woods' covering conjecture has been disproved in dimensions $n \geq 24$. However, it is not clear that the failure of Woods' conjecture in higher dimensions implies the failure of Minkwoski's conjecture. It would be interesting to understand whether Minkowski conjecture holds in some dimensions where Woods' conjecture fails. In addition, the question of which subclasses of lattices satisfy the Woods' conjecture either in dimensions where the question is open or for those where the general conjecture fails is intresting. 
\medskip

\noindent
{\bf Classification of perfect and eutactic lattices} is known only in low dimensions. In fact, there are not even known asymptotic formulas for the number of perfect or eutactic similarity classes of lattices in growing dimensions: the upper and lower bounds known for perfect similarity classes are of different orders of magnitude as functions of the dimension. An important avenue for future research would be to obtain stronger general bounds with a view toward asymptotic formulas.
\medskip

\noindent
{\bf Zeta-function of well-rounded sublattices} in a fixed planar lattice was studied by several authors and its behavior is generally understood. However, there seem to be no analogous results in higher dimensions. Studying the analytic property of this function in higher dimensions would provide insight into the quantitative distribution properties of well-rounded sublattices and their dependence on the arithmetic structure of the ambient lattice.
\medskip

\noindent
{\bf Algebraic constructions of well-rounded lattices} have received some attention in the recent years. In particular, ideal well-rounded lattices from quadratic number fields are fairly well understood. On the other hand, there are few results for number fields of higher degree. Also, constructions of well-rounded lattices from more general free $\zed$-modules in number fields deserve more attention as they can often display interesting properties, such as large automorphism groups. Further investigation of tame lattices is also an interesting related project.
\medskip

\noindent
{\bf Further constructions of lattices with special geometric properties}, such as well-roundedness, stability, eutaxy, and perfection coming from function fields, graph theory, theory of tight frames and other areas of mathematics are of great interest. Connections between lattices and spherical designs is a topic of research that also naturally falls here. Continuing investigations in these directions is certainly worthwhile.
\medskip

\noindent
{\bf Studying the local and global extrema of theta functions and finding close-to-optimal constructions} is a very natural question in mathematics and, as discussed in this survey, also motivated by applications both in lattice-based cryptography and physical layer security via the variational distance. The question of finding the precise minimum and the lattice(s) achieving it is generally very hard, and the answer is known only in a few small dimensions.  Hence, any new insight toward this goal will be valuable. 
\medskip

\noindent
{\bf The equivalence of variants of the LWE problem} such as RLWE and PLWE has been studied in the literature for some classes of number fields. These equivalence results enable the construction of cryptographic  schemes with improved efficiency while maintaining strong security guarantees. Extending these results to broader classes of number fields would provide more options for designing secure and efficient protocols. This will be particularly important if the currently standardized lattice-based schemes (e.g. Kyber) would render themselves vulnerable to fatal attacks. 

\medskip 

\noindent 
{\bf Analyzing the hardness of approximate versions of worst-case hard problems} such as SVP and SIVP for varying approximation factors remains an important research direction. In particular, understanding the hardness of approximate SVP or SIVP over structured lattices may either strengthen the security guarantees of the newly standardized schemes or reveal potential weaknesses in their underlying hardness assumptions.
\medskip 

\noindent
{\bf Cryptanalysis of variants of LWE} has also been extensively studied in the literature. This line of research helps maintain robust security by identifying vulnerable instances before they are exploited by adversaries. Exploring the algebraic structures of the underlying schemes may reveal additional weak instances of these problems.
\medskip 

\noindent
{\bf Constructions of explicit lattice coset codes for wireless communications} under varying channel conditions and dimensions. In particular, it would be interesting to see which further lattice properties (in addition to the density, flatness factor, and product distance) may contribute toward high performance. Moreover, constructing WR lattices from cyclic division algebras and studying their properties in terms of multi-antenna communications would be useful as existing constructions thus far are scarce, especially beyond quaternion algebras.  

\section*{Acknowledgments} We would like to thank Dr. Ragnar Freij-Hollanti and Prof. Russell Lai for useful discussions and helpful comments on the original manuscript.

\bibliographystyle{abbrv}
\bibliography{Survey_main}

\end{document}